\documentclass{tran-l}

 \usepackage[centertags]{amsmath}
 \usepackage{amscd}
 \usepackage{amsfonts}
 \usepackage{amssymb}
 \usepackage{amsthm}
 \usepackage{newlfont}
 \usepackage{graphicx}
 \usepackage[active]{srcltx}
\newcommand{\Real}{\mathbb{R}}

\newcommand{\lc}{\mathcal{L}}
\newcommand{\g}{\mathcal{G}}
\newcommand{\Co}{\operatorname{Cone}}

\theoremstyle{plain}
\newtheorem{thm}{Theorem}[section]
\newtheorem{cor}[thm]{Corollary}
\newtheorem{lem}[thm]{Lemma}
\newtheorem{prop}[thm]{Proposition}
\theoremstyle{definition}
\newtheorem{defn}{Definition}[section]
\theoremstyle{remark}
\newtheorem{rem}{Remark}[section]
\theoremstyle{remark}
\newtheorem{exmp}{Example}[section]
\numberwithin{equation}{section}

\usepackage[cmtip, matrix, arrow]{xy}

\begin{document}

\title{Relative Gerbes}

\author{ Zohreh Shahbazi }

\address{Department of Mathematics, University of Toronto, Toronto, Ontario,
 Canada}

\email{zohreh@math.utoronto.ca}




\dedicatory{}



\begin{abstract}
This paper introduces the notion of ``relative
 gerbes'' for smooth maps of manifolds, and discusses their
 differential geometry. The equivalence classes of relative gerbes
 are further classified by the relative integral cohomology in degree
 three.
\end{abstract}

\maketitle

\section{Introduction}

 \indent Giraud~\cite{Gr} first introduced the concept of gerbes in
the early 1970s to study non-Abelian second cohomology. Later,
Brylinski~\cite{MR94b:57030} defined gerbes as sheaves of groupoids
with certain axioms, and discussed their differential geometry. He
proved that the group of equivalence classes of gerbes gives a
geometric realization of integral three cohomology classes on
manifolds. Through a more elementary approach, Chatterjee and
Hitchin \cite{D,MR2003f:53086} introduced gerbes in terms of
transition line bundles for a given cover of the manifold. From this
point of view, a gerbe is a one-degree-up generalization of a line
bundle, where the line bundle is presented by transition maps. A
notable example of a gerbe arises as the obstruction for the
existence of a lift of a principal $G$-bundle to a central extension
of the Lie group. Another example is the associated gerbe
 of an oriented codimension three submanifold of an oriented
 manifold. The third example is what is called ``basic gerbe,'' which corresponds to the
 generator of the degree three integral cohomology of a compact,
 simple and simply connected Lie group. The basic gerbe over G is closely related to
 the basic central extension of the loop group, and it was constructed, from this point of
 view, by Brylinski~\cite{MR94b:57030}. Later,
 Gawedski-Reis~\cite{MR2003m:81222},
 for G=SU(n), and Meinrenken~\cite{EM}, in the general case, gave a finite-dimensional
 construction along with an explicit description of the gerbe
 connection.\\
\indent This paper introduces the notion of \emph{relative
 gerbes} for smooth maps of manifolds, and discusses their
 differential geometry. The equivalence classes of relative gerbes
 are classified by the relative integral cohomology in degree
 three.\\
 \indent The organization of this paper is as follows. In Section
2, the relative (co)homology of a smooth map between two manifolds
is discussed. When the map is inclusion, the singular relative
(co)homology of the map coincides with the singular relative
(co)homology of the pair. Also, for a continuous map of topological
spaces, the relative (co)homology of the map is isomorphic to the
(co)homology of the mapping cone. In Section 3, following the
Chatterjee-Hitchin perspective on gerbes, the notion of
\emph{relative gerbe} is defined for a smooth map $\Phi\in
C^{\infty}(M,N)$ between two manifolds $M$ and $N$ as a gerbe over
the target space together with a quasi-line bundle for the pull-back
gerbe. It is also proven that the group of equivalence classes of
relative gerbes can be characterized by
the integral degree three relative cohomology of the same map.\\
\indent Another objective of this paper is to develop the
differential geometry of relative gerbes. More specifically, in
Section 4, the concepts of relative connection, relative connection
curvature, relative Cheeger-Simons differential character, and
relative holonomy are introduced. As well, it is proven that a given
closed relative 3-form arises as a curvature of some relative gerbe
with connection if and only if the relative 3-form is integral.
Further, it is shown that a relative gerbe with connection for a
smooth map $\Phi: M\to N$ generates a \emph{relative} line bundle
with connection for the corresponding map of loop paces,
$L\Phi:\,LM\to LN$.

\section
{Relative Homology/Cohomology}
\subsection{Algebraic Mapping Cone for Chain Complexes}
\begin{defn}Let $f_{\bullet}:X_\bullet\longrightarrow Y_\bullet$
be a chain map between chain complexes over R where R is a
commutative ring. The algebraic mapping cone of
$f$~\cite{MR2002f:55001} is defined as a chain complex
$\operatorname{Cone}_{\bullet}(f)$ where
\[\operatorname{Cone}_n(f)=X_{n-1}\oplus Y_{n}\]with the
differential\[\partial(\theta,\eta)=(\partial\theta,f(\theta)-\partial
\eta).\]Since $\partial^2=0$, we can consider the homology of this
chain complex. Define relative homology of $f_{\bullet}$ as
$$H_n(f):=H_n(\operatorname{Cone}_{\bullet}(f)).$$
\end{defn}
The short exact sequence of chain complexes\[0\rightarrow
Y_{n}\overset{j}\rightarrow
\operatorname{Cone}_n(f)\overset{k}\rightarrow X_{n-1}\rightarrow
0\] where $j(\beta)=(0,\beta)$ and $k(\alpha, \beta)=\alpha$ gives a
long exact sequence in
homology\begin{equation}\label{eq0}\begin{split}\cdots \rightarrow
H_{n}(Y)\overset{j}\rightarrow H_{n}(f)\overset{k}\rightarrow
H_{n-1}(X)\overset{\delta}\rightarrow H_{n-1}(Y) \rightarrow
\cdots\end{split}\end{equation}where $\delta$ is the connecting
homomorphism.
\begin{lem}The connecting homomorphism $\delta$ is given by
$\delta[\gamma]=[f(\gamma)]$ for $\gamma\in
X_{n-1}$.\end{lem}\begin{proof} For $\gamma\in X_{n-1}$, we have
$k(\gamma,0)=\gamma$. The short exact sequence of chain complexes
gives an element $\gamma'\in Y_{n-1}$ such that
$j(\gamma')=\partial(\gamma,0)=(\partial\gamma,f(\gamma))$. $\delta$
is defined by $\delta[\gamma]=[\gamma']$. But, by definition of $j$,
$j(\gamma')=(0,\gamma')$. Therefore $f(\gamma)=\gamma'$. This shows
$\delta[\gamma]=[f(\gamma)]$.\end{proof}\begin{defn}We call a chain
map $f_{\bullet}:X_{\bullet}\rightarrow Y_{\bullet}$ a
quasi-isomorphism if it induces isomorphism in cohomology, i.e.,
$H_{\bullet}(X)\overset{\cong}\rightarrow
H_{\bullet}(Y)$.\end{defn}\begin{cor}
$f_{\bullet}:X_{\bullet}\rightarrow Y_{\bullet}$ is a
quasi-isomorphism if and only if
$H_{\bullet}(f)=0$.\end{cor}\begin{proof} $f$ is a
quasi-isomorphism, if and only if the connecting homomorphism in the
long exact sequence \ref{eq0} is an isomorphism.
\end{proof} \begin{defn}A homotopy operator between two chain
complexes
 $f,g:X_{\bullet}\rightarrow Y_{\bullet}$ is a linear map
 $h:X_{\bullet}\rightarrow Y_{\bullet+1}$
 such that \[h\partial+\partial h=f-g.\quad(\star)\]In that case, $f$ and $g$ are
 called chain homotopic and we denote it by $f\simeq
 g$.\newline \indent Two chain maps $f:X_{\bullet}\rightarrow
 Y_{\bullet}$ and $g:Y_{\bullet}\rightarrow
 X_{\bullet}$ are called homotopy inverse if $g\circ
 f\simeq id_X$ and $f\circ
 g\simeq id_Y$ are both homotopic to the identity. If
  $f:X_{\bullet}\rightarrow Y_{\bullet}$ admits a homotopy inverse, it is
  called a homotopy equivalence. In particular, every homotopy equivalence is a
  quasi-isomorphism. \end{defn}
 \begin{prop}Any homotopy between chain maps $f,g:X_\bullet\rightarrow Y_\bullet$
 induces an isomorphism of chain complexes $\operatorname{Cone}(f)_{\bullet}$ and
 $\operatorname{Cone}(g)_{\bullet}$.
\end{prop}
  \begin{proof}Given a homotopy operator $h$ satisfying $(\star\,)$, define a map
  $F:\operatorname{Cone}_\bullet(f)\rightarrow
  \operatorname{Cone}_\bullet(g)$ by
  \[F(\alpha,\beta)=(\alpha,-h(\alpha)+\beta).\]

 Since\[
 \partial F(\alpha,\beta)=(\partial\alpha,g(\alpha)+\partial
 h(\alpha)+\partial\beta)=(\partial\alpha,f(\alpha)-h\partial(\alpha)+\partial\beta)
 =F\partial(\alpha,\beta),\]
 $F$ is a chain map and its inverse map is $F^{-1}(\alpha,\beta)=(\alpha,h(\alpha)+\beta)$.

  \end{proof}
\begin{lem}\label{1}Let \\*\[\begin{CD}
  0 @>>> X_{\bullet}  @>>>  Y_{\bullet} @>>> Z_{\bullet} @>>> 0  \\
   @.  @VVV  @VVV  @VVV  @.\\
   0 @>>>\widetilde{X}_{\bullet}@>>>\widetilde{Y}_{\bullet}@>>>\widetilde{Z}_{\bullet}
   @>>>0
   \end{CD}\]\\*be a commutative diagram of chain maps with exact
   rows. If two of vertical maps are quasi-isomorphisms, then so is
   the third.\end{lem}
   \begin{proof}The statement follows from the 5-Lemma applied to
   the corresponding diagram in homology, \\*\[\begin{CD}
  \cdots @>>> H_{\bullet}(X)  @>>>  H_{\bullet}(Y) @>>> H_{\bullet}(Z) @>>>
  H_{\bullet-1}(X) @>>> \cdots \\
   @.  @VVV  @VVV  @VVV @VVV @.\\
   \cdots @>>>H_{\bullet}(\widetilde{X})@>>>H_{\bullet}(\widetilde{Y})@>>>H_{\bullet}
   (\widetilde{Z})
   @>>>H_{\bullet-1}(\widetilde{X})@>>>\cdots
   \end{CD}\]\end{proof}
  \begin{prop}\label{2}Suppose that we have the following commutative diagram of chain
  maps,\\*\[\begin{CD}
   X_{\bullet}  @>f_{\bullet}>>     Y_{\bullet}\\
   @V\Phi_{\bullet}VV   @V\Psi_{\bullet}VV\\
   \widetilde{X}_{\bullet}@>\widetilde{f}_{\bullet}>>\widetilde{Y}_{\bullet}
   \end{CD}\]\\*such that $\Phi$ and $\Psi$ are quasi-isomorphisms.
   Then the induced map \[F:\operatorname{Cone}_{\bullet}(f)\rightarrow \operatorname{Cone}_
   {\bullet}(\widetilde{f}),\,\,(\alpha,\beta)\mapsto(\Phi(\alpha),\Psi(\beta))\]
   is a quasi-isomorphism.\end{prop}
   \begin{proof}The map $F$ is a chain map
   since,\begin{equation*}\begin{split}
   \partial F(\alpha,\beta)&=\partial(\Phi(\alpha),\Psi(\beta))\\
   & =(\partial\Phi(\alpha),\widetilde{f}(\Phi(\alpha))-\partial\Psi(\beta))\\
   & =(\Phi(\partial\alpha),\Psi(f(\alpha)-\partial\beta))\\
   & =F(\partial\alpha,f(\alpha)-\partial\beta)\\
   & =F\partial(\alpha,\beta).\end{split}\end{equation*}The chain
   map $F$ fits into a commutative diagram,
   \\*\[\begin{CD}
  0 @>>> Y_{\bullet}  @>>>  \Co_{\bullet}(f) @>>> X_{\bullet-1} @>>> 0  \\
   @.  @VVV  @VVV  @VVV  @.\\
   0 @>>>\widetilde{Y}_{\bullet}@>>>\Co_{\bullet}(\widetilde{f})@>>>\widetilde{X}_{\bullet-1}
   @>>>0
   \end{CD}\]\\*Since $\Phi$ and $\Psi$ are quasi-isomorphisms, so is
   $F$ by Lemma \ref{1}.
   \end{proof}
\begin{prop}\label{3}For any chain map $f_\bullet:X_{\bullet}\rightarrow
Y_{\bullet}$,
there is a long exact sequence
\[\cdots\rightarrow H_{n-1}(kerf)\overset{j}\rightarrow H_n(f)\overset{k}
\rightarrow H_n(coker f)\overset{\delta}\rightarrow
H_{n-2}(kerf)\rightarrow H_{n-1}(f)\rightarrow\cdots\]where $j$, $k$
and the connecting homomorphism $\delta$ are defined
by\begin{eqnarray*}j[\theta]
&=& [(\theta,0)] \\ k[(\theta,\eta)] &=& [\eta\, mod\, f(X)] \\
\delta[(\eta\, mod\, f(X))] &=& [\partial\theta]\in
H_{n-2}(kerf).\end{eqnarray*} Here, $\eta\in Y_n$ and
$\partial\eta=f(\theta)$ for some $\theta\in X_{n-1}$. In
particular, if $f$ is an injection then $H_n(f)=H_n(cokerf)$, and if
it is onto then $H_n(f)=H_{n-1}(kerf)$.
\end{prop}\begin{proof} Let $\widetilde{f}_{\bullet}:X_{\bullet}\rightarrow
im(f_\bullet)\subseteq Y_\bullet$ be the chain map $f_\bullet$,
viewed as a map into the sub-complex $f_\bullet(X_\bullet)\subseteq
Y_\bullet$. We have the following short exact sequence
\[0\rightarrow \operatorname{Cone}_n(\widetilde{f})\overset{i}\rightarrow \operatorname{Cone}_n(f)\overset{k}
\rightarrow coker(f_n)\rightarrow 0\] where $k$ is as above and $i$
is the inclusion map. Therefore, there exists a long exact sequence
\begin{equation}
\begin{split}\cdots\rightarrow
H_n(\widetilde{f})\overset{i}\rightarrow
H_n(f)\overset{k}\rightarrow H_n(cokerf) \rightarrow
H_{n-1}(\widetilde{f})\rightarrow\cdots.\end{split}\end{equation}
Let $\widetilde{f}'_\bullet:X_\bullet/kerf_\bullet\rightarrow
im(f_\bullet)$ be the map induced by $f$. Notice that since
$\widetilde{f}'$ is an isomorphism, therefore
$H_\bullet(\widetilde{f}')=0$. By using the long exact sequence
corresponding to the short exact sequence\[0\rightarrow
kerf_{\bullet-1}\overset{\widetilde{j}}\rightarrow
\operatorname{Cone}_\bullet(\widetilde{f})\overset{\pi} \rightarrow
\operatorname{Cone}_\bullet(\widetilde{f}')\rightarrow 0\]where
$\widetilde{j}(\theta)=(\theta,0)$, and
$\pi(\theta,\eta)=(\theta\,mod\,kerf,\eta)$, we see that
$\widetilde{j}$ is a quasi-isomorphism. Since
$j=i\circ\widetilde{j}$, we obtain the long exact sequence
\[\cdots\rightarrow H_{n-1}(kerf)\overset{j}\rightarrow H_n(f)\overset{k}\rightarrow
H_n(cokerf)\rightarrow H_{n-2}(kerf)\rightarrow\cdots.\]To find
connecting homomorphism, assume $[\eta$ mod $f(X)]\in H_n(cokerf)$
for $\eta\in Y_n$. Then $\partial\eta\in f(X)$, i.e.,
$\partial\eta=f(\theta)$ for some $\theta\in X_{n-1}$. Since
\[f(\partial\theta)=\partial f(\theta)=\partial
\partial\eta=0\]then $\partial\theta\in ker(f)$. Also
$k(\theta,\eta)= \eta$ mod $f(X)$ and
$j(\theta)=i\circ\widetilde{j}(\partial\theta)=i(\partial\theta,0)=(\partial\theta,0)=\partial(\theta,\eta)$.
Thus, we have  \[\delta[(\eta\, mod\, f(X))]=[\partial\theta]\in
H_{n-2}(kerf).\]
\end{proof}
\goodbreak
\subsection{Algebraic Mapping Cone for Co-chain Complexes}If
$f^\bullet:X^\bullet\rightarrow Y^\bullet$ is a co-chain map between
co-chain complexes, the algebraic mapping cone of $f$ is defined as
a co-chain complex $\operatorname{Cone}^{\bullet}(f)$ where
\[\operatorname{Cone}^n(f)=Y^{n-1}\oplus X^n\]with the differential
\[d(\alpha,\beta)=(f(\beta)-d\alpha
,d\beta)\]Since $d^2=0$, we can consider the cohomology of this
co-chain complex. Define relative cohomology of $f^{\bullet}$ as
$$H^n(f):=H^n(\operatorname{Cone}^{\bullet}(f)).$$\begin{rem}
Any cochain complex $(X^{\bullet},d)$ may be viewed as a chain
complex $(\widetilde{X}_{\bullet},\partial)$ where
$\widetilde{X}_n=X^{-n}$ and $\partial_n=d^{-n}\quad(n\in
\mathbb{Z})$. This correspondence takes cochain maps
$f^\bullet:X^{\bullet}\rightarrow Y^{\bullet}$ into chain maps
$\widetilde{f}_\bullet:\widetilde{X}_{\bullet}\rightarrow
\widetilde{Y}_{\bullet}$ where $\widetilde{f}_n=f^{-n}$, and
identifies $\Co(\widetilde{f})$ and $\widetilde{\Co(f)}$ up to a
degree shift:\[\begin{array}{c}
   \widetilde{\Co(f)}_n=\Co(f)^{-n}=Y^{-n-1}\oplus
X^{-n}\\
  \Co(\widetilde{f})_n=\widetilde{X}_{n-1}\oplus\widetilde{Y}_n=X^{-n+1}\oplus
Y^{-n}. \\
\end{array}\]Thus, $\widetilde{\Co(f)}_{n}\cong
\Co(\widetilde{f})_{n+1}.$\newline \indent Using this corresondence,
the results for the mapping cone of chain maps are directly carried
over to cochain maps.\end{rem}
\goodbreak
\subsection{Kronecker Pairing}For a chain complex $X_{\bullet}$, the
dual co-chain complex $(X')^{\bullet}$ is defined by
$(X')^{n}=Hom(X_n,R)$ with the dual differential.\begin{prop} Let
$f_\bullet:X_{\bullet}\rightarrow Y_{\bullet}$ be a map between
chain complexes, and $(f')^{\bullet}:(Y')^{\bullet}\rightarrow
(X')^{\bullet}$ be its dual cochain map . Then the bilinear pairing
\[\Co^n(f')\times \Co_n(f)\rightarrow R\]given by the formula\[\langle(\alpha,
\beta),
(\theta,\eta)\rangle=\langle\alpha,\theta\rangle-\langle\beta,\eta\rangle\]for
$(\alpha,\beta)\in \operatorname{Cone}^n(f')$ and $(\theta,\eta)\in
\operatorname{Cone}_n(f)$ induces a pairing in cohomology/homology
\[H^n(f')\times H_n(f)\rightarrow R.\]
\end{prop}
\begin{proof}
It is enough to show that a cocycle paired with a boundary is zero
and a coboundary paired with a cycle is zero. Let $(\alpha,
\beta)=\partial (\alpha', \beta')$ and $\partial (\theta,\eta)=0$.
Therefore, by definition
\[\alpha=f'\beta'-d\alpha',\, \beta=d
\beta'\]and\[\partial \eta=f(\theta), \,\partial \theta=0 .\]
\begin{equation}
\begin{split}
\langle(\alpha,\beta), (\theta,\eta)\rangle &=\langle
\alpha,\theta\rangle-\langle\beta,\eta\rangle\\&=\langle
f'\beta',\theta\rangle-\langle d \alpha',\theta\rangle-\langle d
\beta',\eta\rangle\\&=\langle f'\beta',\theta\rangle-\langle
\alpha',\partial\theta\rangle-\langle
\beta',\partial\eta\rangle\\&=\langle
f'\beta',\theta\rangle-\langle\beta',f(\theta)\rangle\\&=0.
\end{split}
\end{equation} Similarly we can prove that a co-boundary paired with a cycle is
zero.
\end{proof}
\begin{lem}If $f_{\bullet}:X_{\bullet}\rightarrow Y_{\bullet}$ is a chain map
, and $(f')^{\bullet}:(Y')^{\bullet}\rightarrow (X')^{\bullet}$ be
its dual cochain map, then
$\operatorname{Cone}^{\bullet}(f')=(\operatorname{Cone}_{\bullet}(f))'$.\end{lem}
\begin{proof}
Notice that
$\operatorname{Cone}^{n}(f')=(\operatorname{Cone}_{n}(f))'=(X^{n-1})'\oplus
(Y^n)'$. It follows from definitions that
\[\langle d(\alpha,\beta),(\theta,\eta)\rangle=
\langle(\alpha,\beta),\partial(\theta,\eta)\rangle.\]Therefore
differential of $\operatorname{Cone}^{n}(f')$ is dual of
differential of $\operatorname{Cone}_{n}(f)$.
\end{proof}
\goodbreak
\subsection{Singular, De Rham, \v{C}ech Theory}In this Section, two manifolds $M$ and $N$
 and a map $\Phi\in
 C^{\infty}(M,N)$ are fixed.\\  \textbf{Singular} \textbf{relative}
\textbf{homology}:  Consider the push-forward map
$\Phi_{\ast}:S_q(M, R)\rightarrow S_{q}(N, R)$, where R is a
commutative ring and $S_q(M, R), \, S_q(N, R)$ are the singular
chain complexes of M and N respectively. Singular relative homology
is the homology of the chain complex $\Co_{\bullet}(\Phi_*)$, and is
denoted $H_{\bullet}(\Phi,R)$.
\\ \textbf{Singular} \textbf{relative} \textbf{cohomology}:
Consider the pull-back map $\Phi^{\ast}:S^q(N, R)\rightarrow
S^{q}(M, R)$, where R is a commutative ring, and $S^q(M, R)$ and
$S^q(N, R)$ are the singular co-chain complex of M and N,
respectively. Singular relative cohomology is the cohomology of the
co-chain complex $\Co^{\bullet}(\Phi^*)$, and is denoted
$H^{\bullet}(\Phi,R)$.
\\ \textbf{De} \textbf{Rham} \textbf{relative}
\textbf{cohomology}: For $\Phi\in C^{\infty}(M, N)$,  consider the
pull back-map
\[\Phi^*:\Omega^{q}(N)\rightarrow \Omega^q(M)\] between differential co-chain complexes.
In this paper, the cohomology of $\Co^{\bullet}(\Phi^*)$ is denoted
 as $H_{dR}^\bullet(\Phi)$ and called it ``de Rham relative cohomology
 .'' \\ \textbf{\v{C}ech} \textbf{relative}
\textbf{cohomology}: Let $A$ be a R-module, and
$\mathcal{U}=\{U_{\alpha}\}$ be a good cover of a manifold $M$,
i.e., all the finite intersections are contractible. For any
collection of indices $\alpha_0,\cdots,\alpha_p$ such that
$U_{\alpha_0}\cap\cdots\cap U_{\alpha_p}\neq\emptyset$, let
\[U_{\alpha_0\cdots\alpha_p}=U_{\alpha_0}\cap\cdots\cap
U_{\alpha_p}.\]A \v{C}ech-p-cochain $f\in
\check{C}^p(\mathcal{U},A)$ is a function
\[f=\coprod_{\alpha_0\cdots \alpha_p} f_{\alpha_0\cdots \alpha_p}: \coprod_{\alpha_0\cdots
\alpha_p}U_{\alpha_0\cdots \alpha_p}\rightarrow A\] where
$f_{\alpha_0\cdots\alpha_p}$ is locally constant and anti-symmetric
in indices. The differential is defined by
\[(d f)_{\alpha_0\cdots\alpha_{p+1}}=\sum_{i=0}^{p+1}(-1)^if_{\alpha_0\cdots\hat{\alpha}_i
\cdots\alpha_{p+1}}\]where the ``hat'' sign means that the index has
been omitted. Since $d\circ d=0$, one can define \v{C}ech cohomology
groups with coefficients in $A$
as\[\check{H}^p(M,A):=H^p(\check{C}(\mathcal{U},A)).\] Let
$\mathcal{U}=\{{U}_{i}\}_{i\in I} , \, \mathcal{V}=\{{V}_{j}\}_{j\in
J}$ be good covers of M and N, respectively, such that there exists
a map $r:I\rightarrow J$ with $\Phi({U}_{i})\subseteq {V}_{r(i)}$.
Let $\check{C}^\bullet(M, A)$ and $\check{C}^\bullet(N, A)$ be the
\v{C}ech complexes for given covers, where A is an R-module. Using
the pull-back map $\Phi^{\ast}:\check{C}^\bullet(N, A)\rightarrow
\check{C}^\bullet(M,A)$, the relative \v{C}ech cohomology is defined
as the cohomology of $\Co^{\bullet}(\Phi^*)$. Denote this cohomology
by
 $\check{H}^{\bullet}(\Phi, A)$. \newline \indent Suppose that $\underline{A}$
 is one of the sheaves (~\cite{MR95d:14001},~\cite{MR94b:57030})
 $\underline{\mathbb{Z}}, \underline{\mathbb{R}},
 $\underline{U(1)}, $\underline{\Omega^q}$. Denote the space of
k-cochains of the sheaf $\underline{A}$ on $M$ and $N$,
respectively, as $C^k(M,\underline{A})$ and $C^k(N,\underline{A})$.
Here, the differential is defined as above. Again, we have an
induced map
\[ \Phi^*:C^{k}(N,\underline{A})\rightarrow
C^k(M,\underline{A}).\]Denote the  cohomology of
$\Co^\bullet(\Phi^*)$ as $H^*(\Phi,\underline{A})$.\begin{thm}There
is a canonical
 isomorphism $H_{dR}^{n}(\Phi)\cong
 H^{n}(\Phi,\mathbb{R})$.\end{thm}
 \begin{proof} Let $S_{sm}^{\bullet}(M,\mathbb{R})$ and $S_{sm}^{\bullet}(N,\mathbb{R})$
 be the smooth singular cochain complex of $M$ and $N$, respectively
 ~\cite{MR83i:57016}. Consider the following
  diagram :  \\*\[\begin{CD}
\Omega^{n}(N)  @>\Phi^*>>\Omega^{n}(M)\\
@Vg^{n}VV   @Vf^{n}VV\\
S_{sm}^{n}(N, \mathbb{R})  @>\Phi^*>>S_{sm}^{n}(M,\mathbb{R})
\end{CD}\] \\*where $f^n$ is  defined by
$f^n(\omega):\sigma\mapsto\int_{\Delta_n}{\sigma^*\omega}$, for
$\omega\in \Omega^n(M)$ and $\sigma\in S_{n}^{sm}(M)$ is a smooth
singular n-simplex. $g^{\bullet}$ is defined in a similar fashion.
From these definitions, it is clear that the diagram commutes.
$f^{\bullet}$ and $g^{\bullet}$ are quasi-isomorphisms by de Rham
Theorem ~\cite{MR83i:57016}. Define
$k^{\bullet}:\Omega^{\bullet}(\Phi,\mathbb{R})\rightarrow
S_{sm}^{\bullet}(\Phi,\mathbb{R})$ by $k^{\bullet}(\alpha,
\beta)=(f^{\bullet-1}(\alpha), g^\bullet(\beta))$. One can use
Proposition \ref{2} and deduce that $k^\bullet$ is a
quasi-isomorphism. There is a co-chain map
\[l^{\bullet}:S^{\bullet}(M,\mathbb{R})\rightarrow
S^{\bullet}_{sm}(M,\mathbb{R})\]given by the dual of the inclusion
map in chain level. In (~\cite{MR84k:58001}, p.196) it is shown that
$l^{\bullet}$ is a quasi-isomorphism. Therefore, by using
Proposition \ref{2} again,
$$H^{n}(\Phi,\mathbb{R})\cong
H_{sm}^{n}(\Phi,\mathbb{R}).$$Together, one can have
$H^{\bullet}(\Phi,\mathbb{R})\cong
H_{dR}^{\bullet}(\Phi,\mathbb{R})$.

\end{proof}
 \begin{thm}For $\Phi\in C^{\infty}(M,N)$, there is an
isomorphism $H_{dR}^{q}(\Phi)\cong \check{
H}^{q}(\Phi,\mathbb{R})$.\end{thm}
\begin{proof}Let $\mathcal{U}=\{{U}_{i}\}_{i\in I}$ and
$\mathcal{V}=\{{V}_{j}\}_{j\in J}$ be good covers of M and N
together with a map $r:I\rightarrow J$, such that
$\Phi({U}_{i})\subseteq {V}_{r(i)}$. Define the double complex
$E^{p,q}(M)=\check{C}^p(M,\Omega^q)$ where $
\check{C}^p(\mathcal{U},\Omega^q)$ is the set of $q$-forms
$\omega_{\alpha_0\cdots\alpha_p}\in
\Omega^q(U_{\alpha_0\cdots\alpha_p})$ anti-symmetric in indices with
the differential $d$ defined as before. Let
$E^n(M)=\bigoplus_{p+q=n}E^{p,q}(M)$ be the associated total
complex. The map $\Phi:M\rightarrow N$ induces chain maps
$\Phi^*:E^n(N)\rightarrow E^n(M)$. Let's denote the corresponding
algebraic mapping cone as $E^n(\Phi)$. The inclusion
$\check{C}^n(M,\mathcal{U})\rightarrow E^n(M)$ is a
quasi-isomorphism (~\cite{MR83i:57016}, p. 97). There exists a
similar quasi-isomorphism for $N$, and since inclusion maps commute
with pull-back of $\Phi$, one gets a quasi-isomorphism
$\check{C}^n(\Phi)\rightarrow E^n(\Phi)$. Thus, the following
isomorphism is obtained.
\begin{equation}\label{eq1}\check{H}^n(\Phi,\mathbb{R})\cong
H^n(E(\Phi))\end{equation} The map $\Omega^n(M)\rightarrow
E^{0,n}(M)\subset E^n(M)$, given by restrictions of forms
$\alpha\mapsto \alpha|_{U_i}$, is a quasi-isomorphism
(~\cite{MR83i:57016}, p. 96). Again, these maps commute with pull
back, and hence define a quasi-isomorphism
$\Omega^n(\Phi)\rightarrow E^n(\Phi)$ that means
\begin{equation}\label{eq2}H_{dR}^n(\Phi)\cong
H^n(E(\Phi)).\end{equation} By combining Equation \ref{eq1} and
Equation \ref{eq2}, one obtains
$\check{H}^{\bullet}(\Phi,\mathbb{R})\cong H_{dR}^{\bullet}(\Phi)$
\end{proof}\begin{rem} A modification of this argument, working
instead with the double complex $\check{C}^p(M,S^q)$ given by
collection of $S^q(U_{\alpha_0\cdots\alpha_p})$, gives isomorphism
between \v{C}ech relative cohomology and singular relative
cohomology with integer coefficients, hence
$$\check{H}^q(\Phi,\mathbb{Z})\cong H_S^q(\Phi,\mathbb{Z}).$$
\end{rem}
\goodbreak
\subsection{Topological Definition of Relative Homology}Let
$\Phi:M\rightarrow N$
 be an inclusion map,
then the push-forward map $\Phi_*:S_\bullet(M,R)\rightarrow
S_\bullet(N,R)$ is injection. Proposition \ref{3} shows that
$H_{\bullet}(\Phi)\cong H_\bullet(S(N)/S(M))=H_{\bullet}(N,M;R)$. $
H_{\bullet}(N,M;R)$ is known as relative homology. Obviously, this
is a special case of what the author defined as a relative singular
homology of an arbitrary map $\Phi:M\rightarrow N$.

Given a continuous map $f:X\rightarrow Y$ of topological spaces,
define mapping cylinder \[\operatorname{Cyl}_f=\frac{(X\times
I)\sqcup Y}{(x, 1)\sim f(x)}\]and mapping cone ~\cite{MR2002k:55001}
\[\Co_f=\frac{\operatorname{Cyl}_f}{X\times \{0\}}
.\] Let $\Co(X):= X\times I/X\times \{0\}.$
 There are natural maps
$$i:Y\hookrightarrow \Co_f,\,\,\,j:\Co(X)\rightarrow
\Co_f.$$Note that $j$ is an inclusion only if $f$ is an inclusion.
There is a canonical map,
$$h: S_{n-1}(X)\rightarrow S_{n}(\Co(X))$$
with the property $h\circ
\partial+\partial\circ h=k$, where $h$ is defined by replacing a singular n-simplex
with its cone, and $k:X\hookrightarrow\Co(X)$ is the inclusion map.
Define the map $$l_n:\Co_n(f_*)\rightarrow
S_n(\Co_f),\,\,(x,y)\mapsto j_*(h(x))-i_*(y).$$\begin{thm}
$l_{\bullet}$ is a chain map and a quasi-isomorphism. Thus,
$$H_n(f)\cong H_n(\Co_f).$$\end{thm}\begin{proof}Recall that
$\partial(x,y)=(\partial x,f_*(x)-\partial y)$. Since
\begin{equation}\begin{split}l(\partial(0,y))+\partial
l(0,y)&= l((0,-\partial y))-\partial i_*y \\
&= i_*(\partial y)-\partial i_*y\\
&=0\end{split}\end{equation}and
\begin{equation}\begin{split}l(\partial(x,0))+\partial
l(x,0)&= l((\partial x,f(x)))+\partial j_*h(x)\\
&= j_*h(\partial x)-i_*f(x)+\partial j_*h(x)\\
&=j_*k_*(x)-i_*f(x)\\
&=0\end{split}\end{equation}therefore $\partial l+l\partial=0$.
Consider diagram \\*\[\begin{CD}
  0 @>>> S_n(\operatorname{Cyl}_f)  @>>>  S_n(\Co_f) @>>> S_n(\Co_f,\operatorname{Cyl}_f) @>>> 0  \\
   @.  @AAA  @A l AA  @AAA  @.\\
   0 @>>>S_n(Y)@>>>\Co_n(f_*)@>>>S_{n-1}(X)
   @>>>0
   \end{CD}\]\\*where the first row corresponds to the pair
   $(\Co_f,\operatorname{Cyl}_f)$ and the right vertical arrow
   comes from$$S_{n-1}(X)\rightarrow S_n(\Co(X),X)\underset{exision}\cong
   S_n(\Co_f,\operatorname{Cyl}_f).$$The diagram commutes, and the
   rows are exact. Since the right and left vertical maps are
   quasi-isomorphisms, hence so is the middle map.\end{proof}
\goodbreak
\subsection{An Integrality Criterion} If A and B are R-modules, then
any homomorphism $\kappa:A\rightarrow B $ induces homomorphisms
$\kappa:H^n(\Phi, A)\rightarrow H^n(\Phi, B)$ and $\kappa:H_n(\Phi,
A)\rightarrow H_n(\Phi, B)$. In particular, the injection
$\iota:\mathbb{Z}\rightarrow \mathbb{R}$ induces a homomorphism
\[\iota:H^n(\Phi, \mathbb{Z})\rightarrow H^n(\Phi, \mathbb{R})
.\]A class $[\gamma]\in H^n(\Phi, \mathbb{R})$ is called integral in
case $[\gamma]$ lies in the image of the map $\iota$.
\begin{prop}\label{I}A class $[(\alpha, \beta)]\in H^n(\Phi, \mathbb{R})$
is integral if and only if $\int_{\theta}\alpha-\int_{\eta}\beta \in
\mathbb{Z}$ for all cycles
 $(\theta,\eta)\in \Co_n(\Phi, \mathbb{Z})$
. \end{prop}
\begin{proof}
 Consider the following commutative diagram
 \\*\[\begin{CD}0@>>> 0@>>> H^n(\Phi, \mathbb{R})@>{\cong}>>Hom(H_n(\Phi, \mathbb{R}), \mathbb{R})@>>>0\\
@. @AAA @A{\iota}AA @A{\widetilde{\iota}}AA\\
 0@>>>Ext(H_n(\Phi,\mathbb{Z}))@>>>H^n(\Phi, \mathbb{Z})@>{\tau}>>Hom(H_n(\Phi, \mathbb{Z}), \mathbb{Z})@>>>0\\
 \end{CD}\]\\*where
$H^n(\Phi,\mathbb{R})\rightarrow
Hom(H_n(\Phi,\mathbb{R}),\mathbb{R})$ and $\tau$ are pairing given
by integral. The map $\widetilde{\iota}$ is inclusion map,
considering the fact that
$$Hom(H_n(\Phi,\mathbb{R}),\mathbb{R})=
Hom(H_n(\Phi,\mathbb{Z}),\mathbb{R}).$$Thus, $[(\alpha, \beta)]\in
H^n(\Phi, \mathbb{R})$ is integral if
$\int_{\theta}\alpha-\int_{\eta}\beta \in \mathbb{Z}$ for all cycles
 $(\theta,\eta)\in \Co_n(\Phi, \mathbb{Z})$.
 \end{proof}
 \goodbreak
\subsection{Bohr-Sommerfeld Condition}Let $(N,\omega)$ be a
symplectic manifold. Recall that an immersion $\Phi:M\rightarrow N$
is isotropic if $\Phi^* \omega=0$. It is called Lagrangian if
furthermore $\dim M=\frac{1}{2}\dim N$. Suppose that
$H_1(N,\mathbb{Z})=0$ and $\omega$ is integral. A Lagrangian
immersion $\Phi:M\rightarrow N$ is said to satisfy the
Bohr-Sommerfeld condition (~\cite{MR85e:58069},~\cite{MR94g:58085})
if for all 1-cycles $\gamma\in S_1(M)$
$$\frac{1}{2\pi}\int_D\omega\in \mathbb{Z}\,\,\,\,\,\,\mbox{where}\,\,\,\partial
D=\Phi(\gamma).$$ Note that since $\omega$ is integral, the above
condition does not depend on the choice of $D$. Also, if
$\omega=d\theta$ is exact ( for example for the cotangent bundles ),
the condition means that
$$\frac{1}{2\pi}\int_{\gamma}\Phi^*\theta\in \mathbb{Z}\,\,\,\,\,\,\mbox{for all
1-cycles $\gamma$.}$$ In terms of relative cohomology, the above
condition means that $(0,\omega)\in \Omega^2(\Phi)$ defines an
integral class in $H_{dR}^2(\Phi)$. The interesting feature of this
situation is that the forms on $M,N$ are fixed, and it defines a
condition on the map $\Phi$.\begin{exmp}Let
$N=\mathbb{R}^2,\,M=S^1,\,\omega=dx\wedge dy,\,\Phi=\mbox{inclusion
map}$. Then, the immersion $\Phi: S^1\hookrightarrow \mathbb{R}^2$
satisfies the Bohr-Sommerfeld condition.\end{exmp}
\section{Geometric Interpretation of Integral Relative Cohomology
Groups}

Let $\Phi\in C^{\infty}(M,N)$ where M and N are manifolds. Let
$U=\{\emph{U}_{i}\}_{i\in I}$, $V=\{\emph{V}_{j}\}_{j\in J}$ be good
covers of $M$ and $N$ respectively such that there exists a map
$r:I\rightarrow J$ with
$\Phi(\emph{U}_{i})\subseteq\emph{V}_{r(i)}$.
\begin{prop}
$H^q(\Phi, \mathbb{Z})\cong H^{q-1}(\Phi,\underline{U(1)}\,)$ for
$q\geq1$.
\end{prop}
\begin{proof}Consider the following long exact sequence,
\[\cdots\rightarrow H^{q-1}(M,\underline{\Real})\rightarrow
H^q(\Phi,\underline{\Real})\rightarrow
H^q(N,\underline{\Real})\rightarrow
H^q(M,\underline{\Real})\rightarrow \cdots.\]Since
$H^{\bullet}(M,\underline{\Real})=0$ and
$H^{\bullet}(N,\underline{\Real})=0$, one can see that
$H^q(\Phi,\underline{\Real})=0$ for $q>0$. By using the long exact
sequence associated to exponential sequence\[0\rightarrow
\mathbb{Z}\rightarrow
\underline{\mathbb{R}}\overset{\exp}\rightarrow
\underline{U(1)}\rightarrow 0,\,\,\,\,\,\,(\star)\]one can deduce
that $H^q(\Phi,\mathbb{Z})\cong H^{q-1}(\Phi,\underline{U(1)})$ for
$q\geq1$.
\end{proof}\goodbreak
\subsection{Geometric Interpretation of $H^1(\Phi,\mathbb{Z})$} Let
$X$ be a manifold. Function $f\in C^{\infty}(X,U(1))$ has global
logarithm if there exists a function $k\in C^{\infty}(X,\mathbb{R})$
such that $f=\exp((2\pi \sqrt{-1}) k)$.

\begin{defn}
The two maps $f,g:X\rightarrow U(1)$ are equivalent if $f/g$ has a
global logarithm.
\end{defn}
The short exact sequence of sheaves $(\star)$ gives an exact
sequence of Abelian groups\[0\rightarrow
H^0(X,\mathbb{Z})\rightarrow
C^{\infty}(X,\mathbb{R})\overset{\exp}\rightarrow
C^{\infty}(X,U(1))\rightarrow H^1(X,\mathbb{Z})\rightarrow 0.\] This
shows that there is a one-to-one correspondence between equivalence
classes and elements of $H^1(X,\mathbb{Z})$. One should look
 for a geometric realization of $H^1(\Phi,\mathbb{Z})$ for a smooth map $\Phi:M\rightarrow N$. Let\[\mathcal{L}:=
 \{(k,f)|\Phi^*f=\exp((2\pi \sqrt{-1}) k)\}\subset C^{\infty}(M,\mathbb{R})\times
 C^{\infty}(N,U(1)).\]$\mathcal{L}$ has a natural group structure. There is a natural group
homomorphism, \[\tau: C^{\infty}(N,\mathbb{R})\rightarrow
 \mathcal{L}\]where $\tau$
 is defined for $l\in C^{\infty}(N,\mathbb{R})$ by
 \[\tau(l)=(\Phi^*l,\exp((2\pi \sqrt{-1}) l)).\]
 \begin{defn}$(k,f),(k',f')\in \mathcal{L}$ are equivalent if
 $f/f'=\exp((2\pi \sqrt{-1}) h)$ for
 some function $h\in C^{\infty}(N,\mathbb{R})$ such that\[\Phi^*h=k-k'.\]\end{defn}
 The set of equivalence classes is a group
 $\mathcal{L}/\tau(C^{\infty}(N,\underline{\Real}))$.\begin{thm}
There exists an exact sequence of groups
\[C^{\infty}(N,\mathbb{R})\overset{\tau}\rightarrow
 \mathcal{L}\rightarrow H^1(\Phi,\mathbb{Z})\rightarrow 0.\]Thus, $H^1(\Phi,\mathbb{Z})$
 parameterizes equivalence classes of pairs $(k,f)$.\end{thm}\begin{proof} The first step
 is to construct a group homomorphism \[\chi:\mathcal{L}\rightarrow H^1(\Phi,\mathbb{Z}).\]
 Given $(k,f)$, let $l_{j}\in C^{\infty}(V_{j},\mathbb{R})$ be local logarithms
 for $f|_{V_{j}}$, that is $f|_{V_{j}}=\exp((2\pi \sqrt{-1}) l_{j})$. On overlaps,
  $a_{jj'}:=l_{j'}-l_{j} :V_{jj'}\rightarrow \mathbb{Z}$ defines
  a \v{C}ech cocycle in \v{C}$^1(N,\mathbb{Z})$.
  Let \[b_{i}:= \Phi^*l_{r(i)}-k|_{U_{i}}: U_{i}
  \rightarrow \mathbb{Z}.\]Since $b_i'-b_i=\Phi^*a_{r(i)r(i')}$, so that
  $(b,a)$ defines a \v{C}ech cocycle in \v{C}$^1(\Phi,\mathbb{Z})$.
  Given another choice of local logarithms $\widetilde{l}_j$, the
  \v{C}ech cocycle changes to
  \[\widetilde{b}_i=b_i+\Phi^*c_{r(i)},\quad
  \widetilde{a}_{jj'}=a_{jj'}+c_{j'}-c_j\]where
  $c_j=\widetilde{l}_j-l_j: V_j\rightarrow \mathbb{Z}$. Thus,
  $(\, \widetilde{b},\widetilde{a})=(\,b,a)+ d(0,c)$, and
  $\chi(k,f):=[(b,a)]\in H^1(\Phi,\mathbb{Z})$ is well-defined.
  Similarly, if $(b,a)=d(0,c)$ then the new local logarithms
  $\widetilde{l}_j=l_j-c_j$ satisfy $\widetilde{a}_{jj'}=0$, which
  means that $\widetilde{l}_j$ patches to a global logarithm
  $\widetilde{l}$. $b_{i}=\Phi^*c_{r(i)}$ implies that
$k|_{U_{i}}=\Phi^*\widetilde{l}_{r(i)}$, which
  means $k=\Phi^*\widetilde{l}$. This shows that the kernel of $\chi$ consists of $(k,f)$
  such that there exists $l\in C^{\infty}(N,\mathbb{R})$ with $f=\exp((2\pi
\sqrt{-1}) l)$ and $k=\Phi^*l$,
  i.e., $ker(\chi)=im(\tau)$.\newline \indent
  Finally, it is shown below that $\chi$ is surjective. Suppose that
  $(b,a)\in$ \v{C}$^1(\Phi,\mathbb{Z})$ is a cocycle. Then
  \begin{equation}a_{j'j''}-a_{jj''}+a_{jj'}=0\end{equation}
  \begin{equation}\Phi^*a_{r(i)r(i')}=b_{i'}-b_{i}\end{equation}
  Choose a partion of unity $\sum_{j\in J}h_j=1$ subordinate to
  the open cover $V=\{\emph{V}_{j}\}_{j\in J}$. Define $f_j\in C^{\infty}(V_j,U(1))$
  by\[f_{j}=\exp(2\pi \sqrt{-1}\sum_{p\in J}
  a_{{j}{p}}h_p).\]By applying (3.1) on $V_{j}\cap V_{j'}$ one has
  \begin{eqnarray*}f_{j}f_{j'}^{-1}&=&\exp(2\pi \sqrt{-1}\sum_{p\in J}
  a_{{j}{p}}h_p)\exp(-2\pi \sqrt{-1}\sum_{p\in J}
  a_{j' p}h_p)\\
  &=&\exp(2\pi \sqrt{-1}\sum_{p\in J}
  a_{{jj'}}h_p)\\
  &=&1.\end{eqnarray*} Hence $f_i$ defines a map $f\in C^{\infty}(N,U(1))$
  such that $f_{|V_{j}}=f_j$. Define $k_{i}\in C^{\infty}(U_{i},\mathbb{R})$ by
  \begin{equation}k_{i}=
  \sum_{p\in J}( \Phi^*a_{r(i)p}+b_{i})\Phi^*h_p.\end{equation}
  Since $b_{i}\in \mathbb{Z}$
  , $\exp((2\pi \sqrt{-1}) k_{i})=\Phi^*f|_{U_{i}}$. One can check that on overlaps
  $U_{i}\cap U_{i'}$, $k_{i}-k_{i'}=0$,
 so that $\{k_{i}\}$ defines a global function $k\in C^{\infty}(M,\Real)$ with
 $\Phi^*f=\exp((2\pi \sqrt{-1}) k)$. Indeed, by applying (2.1.1)
  and (2.1.2) on $U_{i}\cap U_{i'}$ one can obtain
  \[\sum_{p\in J}(\Phi^*a_{r(i) p}+\Phi^*a_{pr(i')}+b_{i}-b_{i'})\Phi^*h_p=
  \sum_{p\in J}(\Phi^*a_{r(i)r(i')}+b_{i}-b_{i'})\Phi^*h_p=0.\]
  By construction $\chi(k,f)=[(b,a)]$, which shows $\chi$ is surjective.\end{proof}
  \begin{rem}Any $(k,f)\in \mathcal{L}$ defines a $U(1)$-valued
  function on the mapping cone, $\Co_\phi=N\cup_{\Phi}\Co(M)$,
  given by $f$ on $N$ and by $\exp((2\pi \sqrt{-1}) t k)$ on $\Co(M)$.
  Here,
  $t\in I$ is the cone parameter. Hence, one obtains a map $$\mathcal{L}\rightarrow
  H^1(\Co_{\Phi},\mathbb{Z})\cong H^1(\Phi,\mathbb{Z}).$$ This
  gives an alternative way of proving the Theorem
  3.2.\end{rem}\goodbreak\goodbreak
\subsection{Geometric Interpretation of $H^2(\Phi,\mathbb{Z})$}
 Denote the group of Hermitian line bundles over $M$ with $Pic(M)$ and the subgroup
 of Hermitian line bundles over $M$ which admits a unitary section with $Pic_0(M)$.
 Recall that there is an exact sequence of Abelian groups
 \[0\rightarrow Pic_0(M)\hookrightarrow
 Pic(M)\overset{\delta} \rightarrow H^2(M,\mathbb{Z})\rightarrow 0\]
 defined as follows. For the line bundle $L$ over $M$ with transition maps
 $c_{ii'}\in C^{\infty}(U_{ii'},U(1))$ over
 good cover $\{U_i\}_{i\in I}$ for $M$, $\delta(L)$ is the cohomology class
 of the 2-cocycle $\{a_{ii'i''}: U_{ii'i''}\rightarrow \mathbb{Z}$
 given as
 \[a_{ii'i''}:=\big(\frac{1}{2\pi \sqrt{-1}}(\log c_{i'i''}-\log c_{ii''}+
 \log c_{ii'})\big)\in\mathbb{Z}.\]
Thus, one can say two Hermitian line bundles $L_1$ and $L_2$ over
$M$ are equivalent if and only if $L_1L_2^{-1}$ admits a unitary
section. The exact sequence shows that $H^2(M,\mathbb{Z})$
parameterizes the equivalence classes of line
bundles~\cite{MR45:3638}. The class $\delta(L):=c_1(L)$ is called
the first Chern class of $L$. Similarly, for a smooth map
$\Phi:M\rightarrow N$ one should look for a geometric realization of
$H^2(\Phi,\mathbb{Z})$.

 \begin{defn}
 Suppose that $\Phi\in C^{\infty}(M,N)$ and $L_1,L_{2}$ are two Hermitian line bundles
 over N, and $\sigma_{1},\sigma_{2}$ are unitary sections of
 $\Phi^{*}L_1,\Phi^{*}L_{2}$. Then,
 $(L_1,\sigma_1)$ is equivalent to
 $(L_{2},\sigma_{2})$ if $L_1L_2^{-1}$ admits a unitary section $\tau$
 , and there is a map
 $f\in C^{\infty}(M,\mathbb{R})$ such that
 $(\Phi^{*}\tau)/\sigma_1\sigma_2^{-1}=\exp((2\pi \sqrt{-1}) f)$.
 \end{defn}

  This defines an equivalence relation among $(\sigma, L)$, where
  $L$ is a Hermitian line bundle over N and $\sigma$ is a
  unitary section of $\Phi^{*}L$. \begin{defn} A relative line bundle for
  $\Phi\in C^{\infty}(M,N)$ is a pair $(\sigma, L)$, where $L$ is a Hermitian
  line bundle over $N$ and $\sigma$ is an unitary section for $\Phi^*L$. Define the group
of relative line bundles \begin{eqnarray*}
  Pic(\Phi)&=&\{(\sigma, L)|L\in Pic(N),
 \sigma \,\mbox{a unitary section of $\Phi^*L$}\}\\
 \mbox{and a
 subgroup of it}\\
Pic_0(\Phi)&=&\{(\sigma, L)\in Pic(\Phi)|\mbox{$\exists$ a unitary
  section}\, \tau\, \mbox{of L and $k$}\in
C^{\infty}(M,\mathbb{R})\\ & &\mbox{ with
}\Phi^*\tau/\sigma=\exp((2\pi \sqrt{-1})
 k)\}.\end{eqnarray*}\end{defn}
 \begin{exmp}Let $(N,\omega)$ be a compact
 symplectic manifold of dimension 2n, and let $L\rightarrow N$ be a line
 bundle with connection $\nabla$ whose curvature is $\omega$,
 i.e., $L$ is a pre-quantum line bundle with connection. A
 Lagrangian submanifold $M$ satisfies the Bohr-Sommerfeld condition if
 there exists a global non-vanishing covariant constant(=flat)
 section $\sigma_M$ of $\Phi^*L$, where $\Phi:M\rightarrow N$ is
 inclusion map (1.6,~\cite{MR94g:58085}). For any Lagrangian submanifold $M$,
 $(\sigma_M, L)\in Pic(\Phi)$.
 \end{exmp}
  \begin{thm}
  There is a short exact sequence of Abelian groups \[0\rightarrow Pic_0(\Phi)
  \rightarrow Pic(\Phi)\rightarrow H^2(\Phi,\mathbb{Z})\rightarrow 0.\]
  Thus, $H^2(\Phi,\mathbb{Z})$ parameterizes
  the set of equivalence classes of pairs $(\sigma, L)$.
  \end{thm}

\begin{proof}
One can identify $H^2(\Phi,\mathbb{Z})$ with
$H^1(\,\,\Phi\,,\,\underline{U(1)}\,\,)$ by Proposition 3.1. Let
$(\sigma, L)\in Pic(\Phi)$. Let $\{V_{j}\}_{j\in J}$ be a good cover
of N and $\{U_i\}_{i\in I}$ be a good cover of $M$ such that there
exists a map $r:I\rightarrow J$ with $\Phi(U_{i})\subseteq
V_{r(i)}$. Choose unitary sections $\sigma_j$ of $L|_{V_{j}}$. The
corresponding transition functions for $L$ are\[g_{jj'}\in
C^{\infty}(V_{jj'},U(1))\quad j,j'\in J,\quad
g_{jj'}\sigma_{j}={\sigma}_{j'}\quad on \quad V_{jj'}\] Define
$f_i=\Phi^*(\sigma_{r(i)})/\sigma$ on $U_{i}$. Then
\begin{equation}
\begin{split}f_if_{i'}^{-1}&=(\Phi^*(\sigma_{r(i)})/\sigma).
(\Phi^*(\sigma_{r(i')})/\sigma)^{-1}\\&=\Phi^*g_{r(i)r(i')}.\end{split}
\end{equation}
Since \[(\delta g)_{r(i)r(i')r(i'')}=1,\]then
$(f_{i},g_{r(i)r(i')})$ is a cocycle in \v{C}$^1(\Phi,\mathbb{Z})$.
If one changes local sections $\sigma_j$, $j\in J$, then
$(f_i,g_{r(i)r(i')})$ will shift by a co-boundary. Define
$$\chi:Pic(\Phi)\longrightarrow
H^{1}(\,\,\Phi\,,\,\underline{U(1)}\,\,),\qquad (\sigma,
L)\mapsto[(f,g)].$$\\ \indent To find the kernel of $\chi$, suppose
that $(f,g)=\delta(t,c)$. Thus, $g=\delta c$ and
$f=\phi^{*}(c)\exp(2 \pi i h)^{-1}$, where $h$ is the global
logarithm of $t$. Define local section $\tau_j:=\sigma_j/c_j$ on
$V_j$. Since on $V_{jj'}$
\[\sigma_j/c_j=\sigma_{j'}/c_{j'},\] then we obtain a global
section $\tau$. On the other hand,
\[\Phi^*\sigma_{r(i)}/\sigma=f_i=\phi^{*}c_{r(i)}\exp((2\pi \sqrt{-1}) h)^{-1}.\]
Therefore, $\phi^*\tau/\sigma=\exp((2\pi \sqrt{-1}) h)^{-1}.$ This
exactly shows that the kernel of $\chi$ is $Pic_0(\Phi)$.\newline
\indent Next, it will be shown that $\chi$ is onto. Let
$(f_i,g_{jj'})\in C^1(\Phi,\underline{U(1)}\,\,)$ be a cocycle. Pick
a line bundle $L$ over $N$ with $g_{jj'}$ corresponding to local
sections $\sigma_j$. $\Phi^*\sigma_{r(i)}/f_i$ defines local
sections for $\Phi^*L$ over $U_{i}$. On $U_{i}\cap U_{i'}$
\[\Phi^*\sigma_{r(i)}/f_i=\Phi^*\sigma_{r(i')}/f_{i'}\]
which defines a global section $\sigma$ for $\Phi^*L$.
By construction, $\chi(\sigma, L)=[(f,g)]$. This shows that $\chi$
is onto.
 \end{proof}
 \begin{rem}A relative line bundle $(L,\sigma)$ for the map
 $\Phi:M\rightarrow N$ defines a line bundle over the mapping
 cone, $\Co_{\Phi}=N\cup_{\Phi}\Co(M)$. This line bundle is given
 by $L$ on $N\subset \Co_{\Phi}$ and by the trivial line bundle
 on $\Co(M)$. The section $\sigma$ is used to glue these two
 bundles. Hence, one obtains a map $$Pic(\Phi)\rightarrow
 H^2(\Co_{\Phi},\mathbb{Z})\cong H^2(\Phi,\mathbb{Z}).$$\end{rem}

 \goodbreak
 \subsection{Gerbes}The main references for this Section are
 ~\cite{MR2003f:53086},~\cite{D} and ~\cite{1}.\newline \indent
Let $\mathcal{U}=\{U_i\}_{i\in I}$ be an open cover for a manifold
$M$. It will be convenient to introduce the following notations.
 Suppose that there is a collection of line bundles $L_{i^{(0)},...,i^{(n)}}$ on
 $U_{i^{(0)},...,i^{(n)}}$.
 Consider the inclusion maps,\[\delta_k:\emph{U}_{i^{(0)},...,i^{(n+1)}}\rightarrow
 \emph{U}_{i^{(0)},...,\widehat{i^{(k)}},...,i^{(n+1)}}\quad (k=0,\cdots,n+1)\]
 and define Hermitian line bundles $(\delta L)_{i^{(0)},...,i^{(n+1)}}$
 over $\emph{U}_{i^{(0)},...,i^{(n+1)}}$ by\[\delta L:=\underset{k=0}
 {\overset{n+1}{\bigotimes}}
 (\delta_k^{*}L)^{(-1)^k}.\]Notice that $\delta(\delta L)$ is canonically trivial.
 If one has a unitary section $\lambda_{i^{(0)},...,i^{(n)}}$ of $L_{i^{(0)},...,i^{(n)}}$
 for each $U_{i^{(0)},...,i^{(n)}}\neq\emptyset$, then one can define
$\delta\lambda$ in a similar fashion.
 Note that $\delta(\delta\lambda)=1$ as a section of trivial line bundle.
\begin{defn}
A gerbe on a manifold $M$ on an open cover
$\mathcal{U}=\{\emph{U}_{i}\}_{i\in I}$ of $M$ is defined by
Hermitian line bundles $L_{ii'}$ on each $\emph{U}_{ii'}$ such that
$L_{ii'}\cong L_{i'i}^{-1}$, and a unitary section $
\theta_{ii'i''}$ of $\delta L$ on $\emph{U}_{ii'i''}$ such that
$\delta\theta=1$ on $\emph{U}_{ii'i''i'''}$. Denote this data as
$\mathcal{G}=(\mathcal{U},L,\theta)$.\end{defn} Denote the set of
all gerbes on $M$ on the open cover
$\mathcal{U}=\{\emph{U}_{i}\}_{i\in I}$ as $Ger(M,\mathcal{U})$.
Recall that an open cover $\mathcal{V}=\{V_{j}\}_{j\in J}$ is a
refinement of open cover $\mathcal{U}=\{\emph{U}_{i}\}_{i\in I}$ if
there is a map $r:J\rightarrow I$ with $V_j\subset U_i$. In this
case, one gets a map
\[Ger(M,\mathcal{U})\hookrightarrow Ger(M,\mathcal{V}).\]Define
the group of gerbes on $M$ as
\[Ger(M)=\underset{\longrightarrow}\lim \,Ger(M,\mathcal{U}).\]
Define the product of two gerbes $\mathcal{G}$ and $\mathcal{G}'$ to
be the gerbe $\mathcal{G}\otimes \mathcal{G}'$ consisting of an open
cover of $M$, $\mathcal{V}=\{V_{i}\}_{i\in I}$ (common refinement of
open covers of $\mathcal{G}$ and $\mathcal{G}'$), line bundles
$L_{ii'}\otimes L'_{ii'}$ on $V_{ii'}$ and unitary sections
$\theta_{ii'i''}\otimes \theta'_{ii'i''}$ of $\delta(L\otimes L')$
on $V_{ii'i''}$.\newline \indent $\mathcal{G}^{-1}$, the dual of a
gerbe $\mathcal{G}$, is defined by dual bundles $L^{-1}_{ii'}$ on
$U_{ii'}$ and sections $\theta^{-1}$ of $\delta(L^{-1})$ over
$U_{ii'i''}$. Therefore, one get a group structure on $Ger(M)$. If
$\Phi:M\rightarrow N$ be a smooth map between two manifolds and
$\mathcal{G}$ be a gerbe on $N$ with open cover
$\mathcal{V}=\{V_{j}\}_{j\in J}$, the pull-back gerbe
$\Phi^*\mathcal{G}$ is simply defined on
$\mathcal{U}=\{U_{i}\}_{i\in I}$ where $\Phi(U_{i})\subset V_{r(i)}$
for a map $r:I\rightarrow J$, line bundles $\Phi^*L_{r(i)r(i')}$ on
$U_{ii'}$, and unitary sections $\theta$ of $\delta(\Phi^*L)$ on
$U_{ii'i''}$. \goodbreak\begin{defn} (Quasi-line Bundle): A
quasi-line bundle for the gerbe $\mathcal{G}$ on a manifold $M$ on
the open cover $\mathcal{U}=\{U_{i}\}_{i\in I}$ is defined
as:\begin{enumerate}\item a Hermitian line bundle $E_{i}$ over each
$U_{i}$;
\item Unitary sections $\psi_{ii'}$ of
\[(\delta E^{-1})_{ii'}\otimes L_{ii'}\]such that
$\delta\psi=\theta$.\end{enumerate} Denote this quasi-line bundle as
$\mathcal{L}=(E,\psi).$
\end{defn} \begin{prop}Any two quasi-line bundles over a given gerbe differ by
a line bundle.\end{prop} \begin{proof} Consider two quasi-line
bundles
 $\mathcal{L}=(E,\psi)$ and $\widetilde{\mathcal{L}}=(\widetilde{E},\widetilde{\psi})$ for the
 gerbe $\g=(\mathcal{U},L,\theta)$.
 $\psi_{ii'}\otimes \widetilde{\psi}_{ii'}^{-1}$ is a unitary section for
 \begin{eqnarray*}E_{i'}\otimes
E_{i}^{-1}\otimes L^{-1}_{ii'}\otimes \widetilde{E}^{-1}_{i'}\otimes
\widetilde{E}_{i}\otimes L_{ii'}&\cong& E_{i'}\otimes
E_{i}^{-1}\otimes \widetilde{E}^{-1}_{i'}\otimes
\widetilde{E}_{i}\\&\cong& E_{i'}\otimes
\widetilde{E}^{-1}_{i'}\otimes E_{i}^{-1}\otimes
\widetilde{E}_{i}.\end{eqnarray*}Therefore, $E\otimes
\widetilde{E}^{-1}$ defines a line bundle over $M$.\end{proof}

 Denote the group of all
gerbes on $M$ related to the open cover
$\mathcal{U}=\{\emph{U}_{i}\}_{i\in I}$ that admits a quasi-line
bundle as $Ger_0(M,\mathcal{U})$. Define
\[Ger_0(M)=\underset{\longrightarrow}\lim
\,Ger_0(M,\mathcal{U}).\]

\goodbreak
\begin{prop}\label{D.D}
There exists a short exact sequence of groups \[0\rightarrow
Ger_0(M)\hookrightarrow Ger(M)\overset{\chi}\rightarrow
H^{3}(M,\mathbb{Z})\rightarrow 0.\]
\end{prop}
\begin{proof}
Identify $H^{3}(M,\mathbb{Z})$ with $H^{2}(M,\underline{U(1)})$.
Consider the gerbe $\g$ on $M$. Refine the cover such that any
$L_{ii'}$ admits unitary sections $\sigma_{ii'}$. Define
\[t:=(\delta\sigma)\theta^{-1}.\]Thus, $\delta t=1$, which means $t$
is a cocycle. Define\[ \chi(\mathcal{G}):=[t].\] Different sections
shift the cocycle by $\delta$\v{ C}$^{1}(M,\underline{U(1)})$, which
shows that $\chi$ is well-defined. Also, $\chi(\mathcal{G}\otimes
\mathcal{G}')=[tt']=\chi(\mathcal{G})\chi(\mathcal{G}')$, which
proves that $\chi$ is a group homomorphism. Next, it will be shown
that the kernel of $\chi$ is $Ger_0(M)$. For $\mathcal{G}\in
Ger_0(M)$, choose a quasi-line bundle $\mathcal{L}=(E,\psi)$. Thus,
$t=\delta(\sigma\psi^{-1})$. Hence, $\chi(\mathcal{G})=[t]=1$.
Conversely, if $[t]=1$, then\[t=\delta t'\] and by defining the new
sections $\sigma'=t'\sigma$ one infers that
$\delta\sigma'=t\delta\sigma=\theta$, which shows that $\mathcal{G}$
admits a quasi-line bundle.\newline \indent Finally, it will be
shown that $\chi$ is onto. If $ \mathcal{U}=\{U_{i}\}_{i\in I}$ is
an open cover of $M$ and $t_{ii'i''}$ is a cocycle
\v{C}$^{2}(M,\underline{U(1)})$, then define a gerbe $\mathcal{G}$
on $M$ by trivial line bundle $L_{ii'}$ on $U_{ii'}$ and unitary
sections $\sigma_{ii'}$ on $U_{ii'}$. Define $\theta = t\delta\sigma
$. Since $\delta t=1$, then $\delta\theta=1$. By construction,
$\chi(\mathcal{G})=[t]$.

\end{proof}\begin{defn}

Let $\mathcal{G}\in Ger(M)$. $\chi(\mathcal{G})\in
H^2(M,\underline{U(1)})\cong H^3(M,\mathbb{Z})$ is called
Dixmier-Douady class of the gerbe $\mathcal{G}$, which is denoted as
D.D.$\g$.
\end{defn} A gerbe admits a quasi-line bundle if and only if its
Dixmier-Douady class is zero by Proposition \ref{D.D}.
\begin{exmp}\label{Principal bundle}Let $G$ be a Lie group, and $1\rightarrow
U(1)\rightarrow \widehat{G} \overset{\kappa}\rightarrow G\rightarrow
1$ be a central extension. Suppose that $\pi: P\rightarrow M$ is a
principal $G$-bundle. A lift of $\pi: P\rightarrow M$ is a principal
$\widehat{G}$-bundle $\widehat{\pi}: \widehat{P}\rightarrow M$
together with a map $q: \widehat{P}\rightarrow P$ such that
$\widehat{\pi}=\pi\circ q$ and the following diagram commutes:
\\\[\begin{CD}
   \widehat{G}\times \widehat{P}  @>>>   \widehat{P}\\
   @V (\kappa, q)VV   @V qVV\\
   G\times P@ >>> P
   \end{CD}\]\\*
In the above diagram the horizontal maps are respective group
actions. Suppose that $\{U_{i}\}_{i\in I}$ is an open cover of $M$
such that $P|_{U_{i}}:=P_i$ has a lift $\widehat{P}_i$. Define
$\widehat{G}$-equivariant Hermitian line bundles as
$$E_i=\widehat{P_i}\times_{U(1)}\mathbb{C}\longrightarrow P\mid_{U_i}.$$
Since $U(1)$ acts by weight 1 on $E_{i}$, it acts by weight 0 on
$E_{i}\otimes E_{i'}^{-1}:=E_{ii'}$ on $U_{ii'}$. Therefore, $G$
acts on $E_{ii'}$, and $E_{ii'}/G$ is a well-defined Hermitian line
bundle, namely $L_{ii'}$. By construction, $\delta L$ is trivial on
$U_{ii'i''}$, therefore one can pick trivial section $\theta$ that
obviously satisfies the relation $\delta\theta=1$. This shows the
obstruction to lifting $P$ to $\widehat{P}$ defines a gerbe
$\g$.\newline \indent If $E_i\rightarrow U_i$ defines a quasi-line
bundle $\mathcal{L}$ for $\g$, then the line bundles
$\widetilde{E_i}:=E_i\otimes \pi^*L_i^{-1}$ patch together to a
global $\widehat{G}$-equivariant line bundle $\widehat{E}\rightarrow
P$, and the unit circle bundle defines a global lift
$\widehat{P}\rightarrow P$. Conversely, if $P$ admits a global lift
$\widehat{P}$ and $\widetilde{P_i}=:\widehat{P}\mid_{U_i}$, then
$L_{ii'}$ is trivial, which shows that the resulting gerbe is a
trivial one.\end{exmp}
\goodbreak
\begin{exmp}Let $N\subset M$ be an oriented codimension 3
submanifold of an n-oriented manifold $M$. The tubular neighborhood
$U_0$ of $N$ has the form $P\times_{SO(3)}\mathbb{R}^3$ where
$P\rightarrow N$ is the frame bundle. Let $U_1=M-N$. Then, $U_0\cap
U_1\cong P\times_{SO(3)}(\mathbb{R}^3-0)$. Over
$(\mathbb{R}^3-0)\cong S^2\times (0,\infty)$, one has degree 2 line
bundle $E$ that is $SO(3)$ equivariant. Thus,
$$L_{01}:=P\times_{SO(3)}E$$ is a line bundle over $U_0\cap U_1$,
which defines the only transition line bundle. Since there is no
triple intersection, this data defines a gerbe over
$M$.\end{exmp}\goodbreak
  \subsection{Geometric Interpretation of $H^3(\Phi,\mathbb{Z})$}
\begin{defn} A relative gerbe for $\Phi\in C^{\infty}(M,N)$ is a
pair $(\mathcal{L},\mathcal{G})$ where $\g$ is a gerbe over $N$ and
$\mathcal{L}$ is a quasi-line bundle for $\Phi^*\g$.\end{defn}
\textbf{Notation}: Let $\Phi\in C^{\infty}(M,N)$. Then
\begin{eqnarray*}Ger(\Phi)&=&\{(\mathcal{L},\mathcal{G})|(\mathcal{L},\mathcal{G})\,
 \mbox{is a relative gerbe for}\, \Phi\in C^{\infty}(M,N).\}\\
Ger_0(\Phi)&=&\{(\mathcal{L},\mathcal{G})\in
Ger(\Phi)|\mathcal{G}\,\mbox{admits a quasi-line
bundle}\,\mathcal{L'}\,\mbox{s.th the line bundle}\,
\mathcal{L}\otimes\Phi^*\mathcal{L'}^{-1}\\& &\mbox{admits a unitary
section}\}\end{eqnarray*} \goodbreak\begin{exmp}Consider a smooth
map $\Phi:M\rightarrow N$ with $\dim M\leq 2$. Let $\g$ be a gerbe
on $N$. Since $\Phi^*\g$ admits a quasi-line bundle say
$\mathcal{L}$, $(\mathcal{L},\g)$ is a relative gerbe.\end{exmp}
\begin{thm}
There exists a short exact sequence of Abelian groups
$$0\rightarrow Ger_0(\Phi)\hookrightarrow Ger(\Phi)\overset{\kappa}\rightarrow
H^{3}(\Phi,\mathbb{Z})\longrightarrow 0$$
\end{thm}
\begin{proof}
One can identify $H^{3}(\Phi,\mathbb{Z})\cong
H^{2}(\Phi,\underline{U(1)})$. Let $\{V_{j}\}_{j\in J}$ be a good
cover of N and $\{U_i\}_{i\in I}$ be a good cover of $M$ such that
there exists a map $r:I\rightarrow J$ with $\Phi(U_{i})\subseteq
V_{r(i)}$. Let $(\mathcal{L},\mathcal{G})\in Ger(\Phi)$. Refine the
gerbe $\mathcal{G}=(\mathcal{U},L,\theta)$ sufficiently such that
all $L_{jj'}$ admit unitary sections $\sigma_{jj'}$. Then, define
$t_{jj'j''}\in$ \v{C}$^{2}(N,U(1))$ by
\[t:=(\delta\sigma)(\theta)^{-1}.\]Since $\delta\theta=1$
and $\delta(\delta\sigma)=1$, then $\delta t=1$. Let
$\mathcal{L}=(E,\psi)$ be a quasi-line bundle for
$\Phi^{*}\mathcal{G}$ with unitary sections $\psi_{ii'}$ for line
bundles $\big((\delta E)_{ii'}\big)^{-1}\otimes
\Phi^{*}L_{r(i)r(i')}$. Define $s_{ii'}\in$\v{C}$^1(M,U(1))$ by
$$s_{ii'}:=(\psi_{ii'})^{-1}\big((\delta\lambda)_{ii'}^{-1}\otimes
\Phi^{*}\sigma_{r(i)r(i')}\big)$$where $\lambda_{i}$ is a unitary
section for $E_{i}$ . Now\begin{equation}\begin{split}
(\Phi^*t^{-1})\delta(\Phi^{*}\sigma)&=\Phi^{*}\theta
\\&=\delta\psi
\\&=(\delta s)^{-1}(\delta\delta\lambda^{-1}\otimes\delta \Phi^{*}\sigma)\\&=
(\delta s)^{-1}\delta
\Phi^{*}\sigma.\end{split}\end{equation}This proves that $\delta
s=\Phi^*t$. Define the map
\[\kappa:Ger(\Phi)\rightarrow
H^{2}(\Phi,\underline{U(1)}), \quad\kappa(\mathcal{L},\mathcal{G})=
[(s,t)].\] It is straightforward to check that this map is
well-defined, i.e., it is independent of the choice of
$\sigma_{jj'}$ and $\lambda_i$. Conversely, given $[(s,t)]\in
H^{2}\big(\Phi,\underline{U(1)}\big)$, one can pick $\mathcal{G}$
such that $\theta=t^{-1}(\delta\sigma)$ and define
$$\psi_{ii'}=s_{ii'}^{-1}\big((\delta\lambda^{-1})_{ii'}\otimes
\Phi^{*}\sigma_{r(i)r(i')}\big).$$ Since $\delta s=\Phi^*t$, then
$\mathcal{L}=(E,\psi)$ defines a quasi-line bundle for
$\Phi^{*}\mathcal{G}$. The construction shows
$\kappa(\mathcal{L},\mathcal{G})=[(s,t)]$. Therefore $\kappa$ is
onto. \\ \indent It is now shown that $\ker(\kappa)=Ger_0(\Phi)$.
Assume $\kappa(\mathcal{L},\mathcal{G})=[(s,t)]$ is a trivial class.
Therefore, there exists $(\rho,\tau)\in$\v{
C}$^{1}\big(\Phi,\underline{U(1)}\big)$ such that
$(s,t)=\delta(\rho,\tau)=(\Phi^{*}\tau(\delta\rho)^{-1},\delta\tau)$.
 $t=\delta\tau$ shows that $\mathcal{G}$ admits a
quasi-line bundle$\mathcal{L'}$. Thus, $\mathcal{L}\otimes
\Phi^*\mathcal{L'}^{-1}$ defines a line bundle over $M$. The first
Chern class of this line bundle is given by the cocycle
$s(\Phi^*\tau)^{-1}$. The condition $s=(\Phi^*\tau)\delta\rho^{-1}$
shows that this cocycle is exact, i.e., the line bundle
$\mathcal{L}\otimes \Phi^*\mathcal{L'}^{-1}$ admits a unitary
section. Thus, $ker(\kappa)\subseteq Ger_0(\Phi)$. Conversely, if
$(\mathcal{L},\mathcal{G})\in Ger_0(\Phi)$ then the above argument,
read in reverse, shows that $(s,t)$ is exact. Hence,
$Ger_0(\Phi)\subseteq ker(\kappa)$.
\end{proof}\goodbreak
\begin{rem}A relative (topological) gerbe $(\mathcal{L},\g)\in Ger(\Phi)$
defines a (topological)gerbe over the mapping cone by ``gluing'' the
trivial gerbe over $\Co(M)$ with the gerbe $\g$ over $N\subset
\Co_{\Phi}$. Here, the line bundles $E_i$ that define $\mathcal{L}$
play the role of transition line bundles. For gluing the gerbes see
~\cite{S}.\end{rem}\begin{exmp}Let $1\rightarrow U(1)\rightarrow
\widehat{G}\rightarrow G\rightarrow 1$ be a central extension of a
Lie group $G$. Suppose $\Phi\in C^{\infty}(M,N)$ and $Q\rightarrow
N$ is a principal $G$-bundle. If $P=\Phi^*Q\rightarrow M$ admit a
lift $\widehat{P}$, then one obtains an element of
$H^3(\Phi,\mathbb{Z})$.\end{exmp}
\begin{exmp}Suppose that $G$ is a compact Lie group. Recall that the
universal bundle $EG\rightarrow BG$ is a (topological) principal
$G$-bundle with the property that any principal $G$-bundle
$P\rightarrow B$ is obtained as the pull-back by some classifying
map $\Phi:B\rightarrow BG$. While the classifying bundle is infinite
dimensional, it can be written as a limit of finite dimensional
bundles $E_nG\rightarrow B_nG$. For instance, if $G=U(k)$, one can
take $E_nG$ the Stiefel manifold of unitary $k$-frames over the
Grassmanian $Gr_{\mathbb{C}}(k,n)$. Furthermore, if $B$ is given,
any $G$-bundle $P\rightarrow B$ is given by a classifying map
$\Phi:B\rightarrow B_nG$ for some fixed, sufficiently large $n$
depending only on $\dim B$~\cite{MR94k:55001}.\\ \indent It can be
shown that $H^3(BG,\mathbb{Z})$ classifies central extension
$1\rightarrow U(1)\rightarrow\widehat{G}\rightarrow G\rightarrow
1$~\cite{MR97j:58157}. For $n$ sufficiently large,
$H^3(B_nG,\mathbb{Z})=H^3(BG,\mathbb{Z})$. Hence,
$H^3(\Phi,\mathbb{Z})$ classifies pairs $(\widehat{G},\widehat{P})$,
where $\widehat{G}$ is a central extension of $G$ by $U(1)$ and
$\widehat{P}$ is a lift of $\Phi^*EG$ to $\widehat{G}$.\end{exmp}
\goodbreak
\section{Differential Geometry of Relative Gerbes}
\subsection{Connections on Line Bundles}
 Let $L$ be a Hermitian line bundle with Hermitian connection  $\nabla$ over a manifold $M$.
 In terms
 of local unitary sections $\sigma_i$ of $L\mid_{U_i}$ and the
 corresponding transition maps $$g_{ii'}:U_{ii'}\rightarrow
 U(1),$$connection 1-forms $A_{i}$ on $U_{i}$ are defined by
 $\nabla\sigma_i=(2\pi \sqrt{-1}) A_i \sigma_i$. On $U_{ii'}$,
 \[(2\pi \sqrt{-1})(A_{i'}-A_{i})=g_{ii'}^{-1}dg_{ii'}.\]
 Hence, the differentials $dA_i$ agree on overlaps.
 The curvature 2-form $F$ is defined by $F|_{U_{i}}=:dA_{i}$.
 The cohomology class of $F$ is independent of the chosen
connection. The cohomology class of $F$ is the image of the Chern
class $c_1(L)\in H^2(M,\mathbb{Z})$ in $H^2(M,\mathbb{R})$. A given
closed 2-form $F\in\Omega^2(M,\mathbb{R})$ arises as a curvature of
some line bundle with connection if and only if $F$ is
integral~\cite{MR94b:57030}.\newline \indent The line bundle with
connection $(L,\nabla)$ is called flat if $F=0$.
 In this case, define the holonomy of $(L,\nabla)$ as follows.
 Assume that the open cover $\{U_i\}_{i\in I}$ is a good cover of $M$.
 Therefore, $A_{i}=df_{i}$ on $U_{i}$,
  where $f_{i}:U_i\rightarrow \mathbb{R}$ is a smooth map on $U_{i}$.
  Then,\[d(2\pi \sqrt{-1}(f_{i'}-f_{i})-\log
  g_{ii'})=0.\]Thus, $$c_{ii'}:=(2\pi \sqrt{-1}(f_{i'}-f_{i})-\log
  g_{ii'})$$ are constants. Since $\log g$ is only defined modulo
  $2\pi \sqrt{-1}\,\mathbb{Z}$,
  so there exists a collection of constants
  $\widetilde{c}_{ii'}:=c_{ii'}\mod\mathbb{Z}$. Different choices of $f_i$, shift
  this cocyle with a coboundary. The 1-cocycle $\widetilde{c}_{ii'}$ represents a
\v{C}ech class in \v{H}$^1(M,U(1))$, which is called the
\emph{holonomy} of the flat line bundle $L$ with connection
$\nabla$.\newline \indent
 Let $L\rightarrow M$ be
a line bundle with connection $\nabla$, and $\gamma:S^1\rightarrow
M$ a smooth curve. The holonomy of $\nabla$ around $\gamma$ is
defined as the holonomy of the line bundle $\gamma^*L$ with flat
connection $\gamma^*\nabla$. \subsection{Connections on Gerbes}
\begin{defn}
Let $\g=(\mathcal{U},L,\theta)$ be a gerbe on a manifold $M$. A
gerbe connection on $\g$ consists of connections $\nabla_{ii'}$ on
line bundles $L_{ii'}$ such that
$\big(\delta\nabla\big)_{ii'i''}\theta_{ii'i''}:=
\big(\nabla_{i'i''}\otimes\nabla_{ii''}^{-1}\otimes\nabla_{ii'}\big)\theta_{ii'i''}=0$,
together with 2-forms $\varpi_{i}\in\Omega^{2}(\emph{U}_{i})$ such
that on $U_{ii'}$, $$(\delta \varpi)_{ii'}=F_{ii'}=\mbox{the
curvature of}\, \nabla_{ii'}.$$ This connection gerbe is denoted as
 a
pair $(\nabla,\varpi)$. \end{defn}Since $F_{ii'}$ is a closed
2-form, the de Rham differential $\kappa\mid_{U_i}:=d\varpi_{i}$
defines a global 3-form $\kappa$, which is called the
\emph{curvature of the gerbe connection}. $[\kappa]\in
H^3(M,\mathbb{R})$ is the image of the Dixmier-Douady class of the
gerbe under the induced map by inclusion $$ \iota:
H^3(M,\mathbb{Z})\rightarrow H^3(M,\mathbb{R}).$$A given closed
3-form $\kappa\in\Omega^2(M,\mathbb{R})$ arises as a curvature of
some gerbe with connection if and only if i$\kappa$ is
integral~\cite{MR2003f:53086}.
\begin{exmp}Suppose that $\pi:P\rightarrow B$ is a principal
$G$-bundle, and $$1\rightarrow U(1)\rightarrow
\widehat{G}\rightarrow G\rightarrow 1$$a central extension. In
Example \ref{Principal bundle}, a gerbe $\g$ is described whose
Dixmier-Douady class is the obstruction to the existence of a lift
$\widehat{\pi}:\widehat{P}\rightarrow B$. Following
Brylinski~\cite{MR94b:57030} (also see ~\cite{MR1956150}) one can
define a connection on this gerbe. Two ingredients are
required:\newline \indent (i) a principal connection $\theta\in
\Omega^1(P,\mathfrak{g})$,\newline \indent (ii) a splitting
$\tau:P\times_{G} \widehat{\mathfrak{g}}\rightarrow B\times
\mathbb{R}$ of the sequence of vector bundles
$$0\rightarrow B\times \mathbb{R}\rightarrow P\times_{G}
\widehat{\mathfrak{g}}\rightarrow P\times_{G}
\mathfrak{g}\rightarrow 0$$ associated to the sequence of Lie
algebras $0\rightarrow\mathbb{R}\rightarrow
\widehat{\mathfrak{g}}\rightarrow \mathfrak{g}\rightarrow 0.$ For a
given lift $\widehat{\pi}:\widehat{P}\rightarrow B$, with
corresponding projection $q:\widehat{P}\rightarrow P$, one say that
a principal connection $\widehat{\theta}\in
\Omega^1(\widehat{P},\widehat{\mathfrak{g}})$ lifts $\theta$ if its
image under $\Omega^1(\widehat{P},\widehat{\mathfrak{g}})\rightarrow
\Omega^1(\widehat{P},\mathfrak{g})$ coincides with $q^*\theta$.
Given such a lift with curvature
$$F^{\widehat{\theta}}=d\widehat{\theta}+\frac{1}{2}[\widehat{\theta},\widehat{\theta}]\in
\Omega^2(\widehat{P},\widehat{\mathfrak{g}})_{\mbox{basic}}=\Omega^2(B,P\times_G
\mathfrak{g}),$$let
$K^{\widehat{\theta}}:=\tau(F^{\widehat{\theta}})\in
\Omega^2(B,\mathbb{R})$ be its ``scalar part.'' Any pair of lifts of
$(P,\theta)$ differs by a line bundle with connection $(L,\nabla^L)$
on $B$. Twisting a given lift $(\widehat{P},\widehat{\theta})$ by
such a line bundle, the scalar part changes by the curvature of the
line bundle~\cite{MR94b:57030}
\begin{equation}\label{equation}K^{\widehat{\theta}}+\frac{1}{2\pi
\sqrt{-1}}curv(\nabla^L).\end{equation}In particular, the exact
3-form $dK^{\widehat{\theta}}\in \Omega^3(B)$ only depends on the
choice of splitting and the connection $\theta$. (It does not depend
on choice of lift.) In general, a global lift $\widehat{P}$ of $P$
does not exist. However, let us choose local lifts
$(\widehat{P_i},\widehat{\theta}_i)$ of $(P\mid_{U_i},\theta)$.
Denote the scalar part of $F^{\widehat{\theta}_i}$ with $\varpi_i\in
\Omega^2(U_i)$, and let $L_{ii'}\rightarrow U_{ii'}$ be the line
bundle with connection $\nabla^{L_{ii'}}$ defined by two lifts
$(\widehat{P_i}\mid_{ U_{ii'}},\widehat{\theta}_i)$ and
$(\widehat{P}_{i'}\mid_{U_{ii'}},\widehat{\theta}_{i'})$. By
Equation
\ref{equation},$$(\delta\varpi)_{ii'}=\frac{1}{2\pi\sqrt{-1}}curv(\nabla^{L_{ii'}}).$$On
the other hand, the connection $\delta\nabla^L$ on $(\delta
L)_{ii'i''}=L_{ii'}L_{ii''}^{-1}L_{ii'}$ is just the trivial
connection on the trivial line bundle. Hence, a gerbe connection is
defined.
\newline \indent A quasi-line bundle $(E,\psi)$ with connection
$\nabla^E$ for this gerbe with connection gives rise to a global
lift $(\widehat{P},\widehat{\theta})$ of $(P,\theta)$, where
$\widehat{P}\mid_{U_i}$ is obtained by twisting $\widehat{P}_i$ by
the line bundle with connection $(E_i,\nabla^{E_i})$. The error
2-form is the scalar part of $F^{\widehat{\theta}}$. (see Definition
4.2)\end{exmp}
\begin{defn}Let $\mathcal{G}$ be a gerbe with connection with a
quasi-line bundle $\mathcal{L}=(E,\psi)$. A connection on a
quasi-line bundle consists of connections $\nabla_i^E$ on line
bundles $E_{i}$ with curvature $F^E_i$ such that
\[(\delta\nabla^E)_{ii'}:=\nabla_{i'}^E\otimes(\nabla_i^E)^{-1}\cong
\nabla_{ii'}.\]Also, the 2-curvatures obey $(\delta
F^E)_{ii'}=F_{ii'}$. Denote this quasi-line bundle with connection
by $(\mathcal{L},\nabla^E)$. Locally defined
 2-forms $\omega\mid_{U_i}=\varpi_{i}-F^E_{i}$ patch together
 to define a global 2-form $\omega$,
 which is called the error 2-form~\cite{1}.\end{defn}
 \begin{rem}\label{qc}The difference
between two quasi-line bundles with connections is a line bundle
with connection, with the curvature equal to the difference of the
error 2-forms.\end{rem}

 Let $\mathcal{G}=(\mathcal{U},L,\theta)$ be a gerbe with connection on $M$.
 Again, assume that $\mathcal{U}$
 is a good cover. Let
 $t\in$\v{C}$^2\big(M,\underline{U(1)}\big)$ be a representative for the Dixmier-Douady
 class of $\g$. Then,
  one canhave a collection of 1-forms $A_{ii'}\in \Omega^1(U_{ii'})$ and 2-forms
  $\varpi_{i}\in \Omega^2(U_{i})$ such that
  \[\kappa|_{U_{i}}=d\varpi_{i}\]\[\delta \varpi=dA\]
  \[(2\pi \sqrt{-1})\delta A=t^{-1}dt.\]If $\kappa=0$, the gerbe is called flat. In this
  case by using Poincar\'e Lemma,
  $\varpi_i=d\mu_i$ on $U_i$ and on $U_{ii'}$,
  $$(\delta\varpi)_{ii'}=d\delta(\mu)_{ii'}=dA_{ii'}.$$Thus, again by Poincar\'e Lemma
  $$A_{ii'}-(\delta\mu)_{ii'}
  =dh_{ii'}.$$ By using $(2\pi \sqrt{-1})\delta A=t^{-1}dt$,
  $$d((2\pi \sqrt{-1})\delta h-\log t)=0.$$
 Therefore, there exists a collection of constants $c_{ii'i''}
  \in$\v{C}$^2(M,\mathbb{R})$. Since $\log$ is defined
  modulo $2\pi \sqrt{-1}\,\mathbb{Z}$, we define
  $$\widetilde{c}_{ii'i''}:=c_{ii'i''}\mod \mathbb{Z}.$$
  The 2-cocycle $\widetilde{c}_{ii'i''}$ represents
  a \v{C}ech class in \v{H}$^2(M,U(1))$, which is called the \emph{holonomy of the flat
  gerbe with connection}.
  Let $\sigma:\Sigma \rightarrow M$ be a smooth map, where
$\Sigma$ is a closed
  surface. The holonomy of $\mathcal{G}$ around $\Sigma$ is defined as the holonomy of
  the pullback gerbe $\sigma^*\mathcal{G}$ with the flat connection
  $\sigma^*(\nabla,\varpi)$(~\cite{MR2003f:53086},~\cite{MR1932333}).
  \subsection{Connections on Relative Gerbes}
  Let $\Phi\in C^{\infty}(M,N)$ and
$U=\{\emph{U}_{i}\}_{i\in I}$, $V=\{\emph{V}_{j}\}_{j\in J}$ are
good covers of $M$ and $N$ respectively such that there exists a map
$r:I\rightarrow J$ with
$\Phi(\emph{U}_{i})\subseteq\emph{V}_{r(i)}$. \begin{defn} A
relative connection on a relative gerbe $(\mathcal{L},\g)$ consists
of gerbe connection $(\nabla,\varpi)$ on $\g$ and a connection
$\nabla^E$ on the quasi-line bundle $\mathcal{L}=(E,\psi)$ for the
$\Phi^*\g$.\end{defn}Consider a relative connection on a relative
gerbe $(\mathcal{L},\g)$. Define the 2-form $\tau$ on $M$ by
$$\tau\mid_{U_i}:=\Phi^*\varpi_{r(i)}-F^E_i.$$ Thus, $(\tau,\kappa)\in \Omega^3(\Phi)$ is
a relative closed 3-form which is called here the \emph{curvature of
the relative connection}.\begin{thm}\label{quantization}A given
closed relative 3-form $(\tau,\kappa)\in \Omega^3(\Phi)$ arises as a
curvature of some relative gerbe with connection if and only if
$(\tau,\kappa)$ is integral.\end{thm}\begin{proof}Let
$(\tau,\kappa)\in \Omega^3(\Phi)$ be an integral relative 3-form. By
Proposition
\ref{I},\begin{equation}\label{This}\int_{\alpha}\kappa-\int_{\beta}\tau\in
\mathbb{Z},\end{equation} where $\alpha\subset N$ is a smooth
3-chain and $\Phi(\beta)=\partial \alpha$, i.e., $(\beta,\alpha)\in
\Co_3(\Phi,\mathbb{Z})$ is a cycle. If $\alpha$ is a cycle then
$(0,\alpha)\in \Co_3(\Phi,\mathbb{Z})$ is a cycle. In this case,
equation \ref{This} shows that for all cycles $\alpha\in
S_3(N,\mathbb{Z})$,
$$\int_{\alpha}\kappa\in \mathbb{Z}.$$ Therefore, one can pick a
gerbe $\g=(\mathcal{V},L,\theta)$ with connection $(\nabla,\varpi)$
over $N$ with curvature 3-form $\kappa$. Denote
$\tau_i:=\tau\mid_{U_i}$. Define $F_i^E\in\Omega^2(U_i)$ by
$$F_i^E=\Phi^*\varpi_{r(i)}-\tau_i.$$ Let
$(\alpha_i,\beta_i)\in\Co_3(\Phi|_{U_i},\mathbb{Z})$ be a cycle.
Then,
\begin{eqnarray*} \int_{\beta_i}F_i^E &=&
\int_{\beta_i}(\Phi^*\varpi_{(r(i)}-\tau_i)\\&=&
\int_{\Phi(\beta_i)}\varpi-\int_{\beta_i}\tau\\&=&\int_{\alpha_i}d\varpi-\int_{\beta_i}\tau\\
&=&\int_{\alpha_i}\kappa-\int_{\beta_i}\tau\in
\mathbb{Z}.\end{eqnarray*}Therefore, one can find a line bundle
$E_i$ with connection over $U_i$ whose curvature is equal to
$F_i^E$. Over $U_{ii'}$, the curvature of two line bundles
$\Phi^*L_{ii'}$ and $E_{i'}\otimes E_i^{-1}$ agrees. Assume that the
open cover $\mathcal{U}=\{U_i\}_{i\in I}$ is a good cover of $M$.
Thus, there is a unitary section $\psi_{ii'}$ for the line bundle
$E_{i}\otimes E_{i'}^{-1}\otimes \Phi^*L_{ii'}$ such that
$\delta\psi=\Phi^*\theta$. Therefore, one obtains a quasi-line
bundle $\mathcal{L}=(E,\psi)$ with connection for $\Phi^*\g$. By
construction the curvature of the relative gerbe $(\lc,\g)$ is
$(\tau,\kappa)$. Conversely, for a given relative gerbe with
connection $(\lc,\g)$ one can have $\int_{\beta_i}F_i^E\in
\mathbb{Z}$ where $\beta_i\subset U_i$ is a 2-cycle which gives
\ref{This}.
\end{proof} Suppose that $\mathcal{G}$ is a gerbe with a flat
connection $(\nabla,\varpi)$
   on $N$ and $\mathcal{L}$
  a quasi-line bundle with connection for $\Phi^*\mathcal{G}$.
  Since $\kappa=0$, as explained in the previous section, there exist 2-cocycles
  $\widetilde{c}_{ii'i''}$ that represent a cohomology class in
  \v{H}$^2(M,U(1))$. Since $\Phi^*\mathcal{G}$
  is trivializable, there is a collection of maps $f_{ii'}$ on $U_{ii'}$
  such that $\delta f=\Phi^*t$ where $j=r(i)$ and $j'=r(i')$. Define $k_{ii'}\in \mathbb{R}$
  as
  \[k_{ii'}=:(2\pi \sqrt{-1})\Phi^*h_{ii'}-\log f_{ii'},\]and \[\widetilde{k}_{ii'}:=
  k_{ii'}\mod\mathbb{Z}.\]Thus, \[\Phi^*\widetilde{c}=\delta \widetilde{k} .\]
  Define the \emph{relative holonomy}
   of the pair $(\,\mathcal{G},\mathcal{L})$ by the relative class
   $[(\widetilde{k},\widetilde{c})]\in H^2(\Phi,U(1))$.

   \begin{defn}
   Let the following diagram be commutative:\\*\[\begin{CD}
   S^1   @>i>>     \Sigma\\
   @VV\psi V   @VV\widetilde{\psi}V\\
   M@>\Phi>>N
   \end{CD}\]\\where $\Sigma$ is a closed surface, i is inclusion map and all
   other maps are smooth. Suppose that $\mathcal{G}$ is a gerbe with connection on $N$, and
   $\Phi^*\mathcal{G}$ admits a quasi-line bundle $\mathcal{L}$ with
   connection. Clearly, $\widetilde{\psi}^*\mathcal{G}$ is a flat gerbe and since
   $i^*\widetilde{\psi}^*\mathcal{G}=\psi^*\Phi^*\mathcal{G}$ then
   $i^*\widetilde{\psi}^*\mathcal{G}$ admits a quasi
   line bundle with connection that is equal to $\psi^*\mathcal{L}$ . Define the
   holonomy of the relative gerbe around the commutative diagram
   as holonomy of the pair $(\psi^*\mathcal{L},\widetilde{\psi}^*\mathcal{G})$.
  \end{defn}
\subsection{Cheeger-Simons Differential Characters}In this section,
a relative version of Cheeger-Simons differential characters is
developed~\cite{MR2031913}~\cite{math.AT/0408333}~\cite{MR1833777}~\cite{MR1976955}.
Denote the smooth singular chain complex on a manifold $M$ as
$S_{\bullet}^{sm}(M)$. Let $Z_{\bullet}^{sm}(M)\subseteq
S_{\bullet}^{sm}(M)$ be the sub-complex of smooth cycles. Recall
that a differential character of degree $k$ on a manifold $M$ is a
homomorphism
$$j:Z_{k-1}^{sm}(M)\rightarrow U(1),$$ such that there is a
closed form $\alpha\in\Omega^k(M)$ with
$$j(\partial x)=\exp\big(2\pi \sqrt{-1}\int_{x}\alpha\big)$$for any
$x\in S_k^{sm}(M)$~\cite{csdc}.

 A connection on a line bundle
defines a differential character of degree 2, where $j$ is the
holonomy map. Similarly, a connection on a gerbe defines a
differential character of degree 3. Specifically, any smooth
$k$-chain $x\in S_k^{sm}(M)$ is realized as a piecewise smooth map
$$\varphi_x:K_x\rightarrow M,$$where $K_x$ is a $k$-dimensional
simplicial complex~\cite{MR2002k:55001}. Then, by definition
$$\int_{K_x}\alpha= \int_x\alpha, \qquad \alpha\in \Omega^k(M).$$
Suppose that $y=\Sigma\epsilon_i\sigma_i\in Z_2^{sm}(M)$, where
$\epsilon_i=\pm 1$. Assume that $\g$ is a gerbe with connection over
$M$. Since $H^3(K_y,\mathbb{Z})=0,$ $\varphi_y^*\g$ admits a
piecewise smooth quasi-line bundle $\mathcal{L}$ with connection.
That is, a quasi-line bundle $\mathcal{L}_i$ for all
$\varphi^*\g\mid_{\Delta^k_{\sigma_i}}$, such that all
$\mathcal{L}_i$ agree on the matching boundary faces. Let $\omega\in
\Omega^2(K_y)$ be the error 2-form and define
$$j(y):=\exp\big(2\pi \sqrt{-1}\int_{K_{y}}\omega\big).$$Any two
quasi-line bundles differ by a line bundle, and hence different
choices for $\mathcal{L}_i$, change $\omega$ by an integral 2-form.
Therefore, $j$ is well-defined. Assume that $y=\partial x$. Since
the components of $K_x$ with empty boundary do not contribute, one
can assume that each component of $K_x$ has non-empty boundary.
Since $H^3(K_x,\mathbb{Z})=0$, choose a quasi-line bundle with
connection for $\varphi^*_x\g$ with error 2-form $\omega$. Let $k$
be the curvature of $\g$. Since $\varphi^*_xk=d\omega$, by stokes'
Theorem
$$\int_{K_x}k=\int_{K_x}d\omega=\int_{\partial
K_x}\omega=\int_{K_y}\omega.$$This shows that $j$ is a differential
character of degree 3.

\begin{defn}Let $\Phi\in C^{\infty}(M,N)$ be a smooth map between
manifolds. A relative differential character of degree $k$ for the
map $\Phi$ is a homomorphism
$$j:Z_{k-1}^{sm}(\Phi)\rightarrow U(1),$$ such that there is a
closed relative form $(\beta,\alpha)\in\Omega^k(\Phi)$ with
$$j(\partial(y,x))=\exp\big(2\pi \sqrt{-1}(\int_{y}\beta-\int_{x}\alpha)\big)$$for any
$(y,x)\in S_k^{sm}(\Phi)$.\end{defn}\begin{thm}A relative connection
on a relative gerbe defines a relative differential character of
degree 3.\end{thm} \begin{proof}Let $\Phi\in C^{\infty}(M,N)$ be a
smooth map between manifolds, and consider a relative gerbe
$(\mathcal{L},\g)$ with connection. Let $(y,x)\in S^{sm}_1(\Phi)$ be
a smooth relative singular cycle, i.e.,
$$\partial y=0$$and $$\Phi_*(y)=\partial x.$$ Suppose that $K_y$ and $K_x$
are the corresponding simplicial complex, and $$\Phi:K_y\rightarrow
K_x$$ is the induced map. Given a relative connection, choose a
quasi-line bundle $\mathcal{L'}$ for $\varphi_x^*\g$, and a unitary
section $\sigma$ of the line bundle
$H:=\varphi_y^*\mathcal{L}\otimes(\Phi^*\mathcal{L}')^{-1}$. Let
$\widetilde{\omega}\in \Omega^2(N)$ be the error 2-form for
$\mathcal{L}'$, and $A\in \Omega^1(M)$ be the connection 1-form for
$H$ with respect to $\sigma$. Define a map $j$ by
$$j(y,x):=\exp\big(2\pi \sqrt{-1}\big(\int_{K_x}\widetilde{\omega}-\int_{K_y}A\big)\big).$$
Choose another quasi-line bundle for $\varphi_x^*\g$. Then, the
difference of error 2-forms is an integral 2-form. Changing the
section $\sigma$ will shift connection 1-form $A$ to $A+A'$, where
$A'$ is an integral 1-form. Thus,
$$j:Z_2^{sm}(\Phi)\rightarrow U(1)$$ is well-defined. Let
$k$ be the curvature 3-form for $\g$, and $\omega$ be the error
2-form for $\mathcal{L}$. Then $(\omega,k)\in \Omega^3(\Phi)$, and
\begin{eqnarray*}j\big(\partial(y,x)\big)&=&j(\partial
y,\Phi_*(y)-\partial x)\\
&=&\exp\big(2\pi \sqrt{-1}\big(\int_{K_{(\Phi_*(y)-\partial
x)}}\widetilde{\omega}-\int_{K_{\partial
y}}A\big)\big)\\&=&\exp\big(2\pi
\sqrt{-1}\big(\int_{K_{\Phi_*(y)}}\widetilde{\omega}
-\int_{K_{\partial x}}\widetilde{\omega}-\int_{K_{\partial
y}}A\big)\big)\\&=&\exp\big(2\pi
\sqrt{-1}\big(\int_{K_{y}}(\Phi^*\widetilde{\omega}-dA)
-\int_{K_x}d\widetilde{\omega}\big)\big)\\
&=&\exp\big(2\pi \sqrt{-1}\big(\int_{K_{y}}\omega
-\int_{K_x}k\big)\big).\end{eqnarray*}Thus, $j$ is a relative
differential character in degree 3.\end{proof}
\subsection{Transgression} For a manifold $M$, denote its loop
space as $LM$. In this Section, a line bundle with connection over
$LM$ is first constructed by transgressing a gerbe with connection
over $M$. A map $\Phi\in C^{\infty}(M,N)$ induces a map $L\Phi\in
C^{\infty}(LM,LN)$. Next, it is proven that a relative gerbe with
connection on $\Phi$ produces a relative line bundle with connection
on $L\Phi$ by transgression. \begin{prop}(Parallel transportation)
Suppose that $\g$ is a gerbe with connection on $M\times [0,1]$ and
$\g_0=\g|_{(M\times \{0\})}$. There is a natural quasi-line bundle
with connection for the gerbe $\pi^*\g_0\otimes\g^{-1}$, where $\pi$
is the projection map
\[\pi:M\times [0,1]\rightarrow M\times \{0\}.\]\end{prop}\begin{proof}
It is obvious that one can obtain a quasi-line bundle with
connection for the gerbe $\pi^*\g_0\otimes\g^{-1}$. Specify a
quasi-line bundle ${\lc}_{\g}$ for the gerbe
$\pi^*\g_0\otimes\g^{-1}$ by the following requirements:\newline
\indent 1.The pull-back $\iota^*{\lc}_{\g}$ is trivial, while
$\iota$ is inclusion map
\[\iota: M\times \{0\}\hookrightarrow M\times [0,1].\]2. Let
$\eta\in\Omega^3(M\times [0,1])$ be the curvature 3-form for
$\pi^*\g_0\otimes\g^{-1}$. Note that $\iota^*\eta=0$. Let
$\chi\in\Omega^2(M\times [0,1])$ be the canonical primitive of
$\eta$ given by transgression. Then, choose a connection on
${\lc}_{\g}$ such that its error 2-form is $\chi$. Any two such
quasi-line bundles differ by a flat line bundle over $M\times
[0,1]$. This line bundle is a trivial line bundle over $M\times
\{0\}$.\end{proof}
\begin{thm}\label{10} A gerbe $\g$ with connection
on $M\times S^1$, induces a line bundle $E_{\g}$ with connection on
$M$. Also, a quasi-line bundle with connection for $\g$ induces a
unitary section of $E_{\g}$.\end{thm}\begin{proof}$M\times
S^1=M\times [0,1]/\sim$, where the equivalence relation is defined
by $(m,0)\sim (m,1)$ for $m\in M$. Therefore,
$\pi^*\g_0\otimes\g^{-1}|_{M\times \{1\}/\sim}$ is a trivial gerbe,
and $\lc_{\g}|_{M\times \{1\}/\sim}$ is a quasi-line bundle with
connection for this trivial gerbe, i.e., a line bundle with
connection $E_{\g}$ for $M$. If one change $\lc_{\g}$ to another
natural quasi-line bundle with connection, the difference between
two quasi-line bundles over $M\times S^1$ is a trivial line bundle.
Thus, the assignment $\g \rightarrow E_{\g}$ is well-defined.
\newline \indent Suppose that the gerbe $\g$ admits a quasi-line bundle $\lc$.
Then, $(\pi^*\lc_0)\otimes(\lc^{-1})$ and $\lc_{\g}$ are two
quasi-line bundles for the gerbe $\pi^*\g_{0}\otimes\g^{-1}$, where
$\lc_0=\lc|_{M\times\{0\}}$. Thus,
$\pi^*\lc_{0}\otimes\lc^{-1}\otimes(\lc_{\g})^{-1}$ defines a line
bundle over $M\times S^1=M\times [0,1]/\sim$. This line bundle over
$M$ defines a map $s:M\rightarrow U(1)$.
$(\pi^*\lc_0)\otimes(\lc^{-1})|_{M\times\{0\}/\sim}$ is the trivial
line bundle $E$. Since $E_{\g}\otimes E^{-1}=s$, $E_{\g}$ admits a
unitary section.\end{proof}
\begin{rem}Let $\g$ be a gerbe with connection on $M$. Consider the
evaluation map \[e:LM\times S^1\rightarrow M.\]Thus, $e^*\g$ induces
a line bundle with connection on $LM$.\end{rem}\begin{thm} For a
given map $\Phi\in C^{\infty}(M,N)$, a relative gerbe with
connection $\g_{\Phi}$ induces a relative line bundle with
connection $E_{L\Phi}$.\end{thm}\begin{proof}The relative gerbe
$\g_{\Phi}$ is a gerbe $\g$ on $N$ together with a quasi-line bundle
with connection $\lc$ for the pull-back gerbe $\Phi^*\g$. The gerbe
$\g$ induces a line bundle with connection $E_{\g}$. Further, the
quasi-line bundle with connection $\lc$ for $\Phi^*\g$ induces a
unitary section $s$ for the line bundle with connection
$(L\Phi)^*E_{\g}$ by Theorem \ref{10}. Thus, the pair $(s,E_{\g})$
defines a relative line bundle with connection
$E_{L\Phi}$.\end{proof}



\subsection*{Acknowledgment}
The results of this paper were obtained during the author's Ph.D.
study at the University of Toronto and were also included in her
thesis dissertation with the same title. The guidance and support
received from Eckhard Meinrenken are deeply acknowledged.

\bibliographystyle{amsplain}
\bibliography{xbib}
\end{document}